\documentclass[a4paper,11pt]{article}
\usepackage{amsmath}
\usepackage{latexsym}
\usepackage{amsfonts}
\usepackage{a4wide}
\usepackage{xspace}
\usepackage{color}
\usepackage{accents}

\usepackage[all]{xy}

\newcommand{\hfl}{\hspace*{\fill}}

\newlength{\dede}     
\newcommand{\tab}[1]{
\settowidth{\dede}{#1}
\hspace*{\dede}}      

\newcommand{\nhp}{\hspace*{-\parindent}}
\newcommand{\noi}{\hspace*{-\parindent}}

\newcommand{\Zh}{\mbox{$\mathbb Z$}} 
\newcommand{\zt}{\mbox{$\se {\Zh}{2}$}}
\newcommand{\Ft}{\mbox{$\se {\mathbb F} 2$}}

\newcommand{\po}[2]{\mbox{$#2^{#1}$}}

\newcommand{\lra}{\mbox{$\longrightarrow$}}  
\newcommand{\com}{\mbox{$\circ$}} 

\newcommand{\se}[2]{\mbox{$ #1_{\raisebox{-3pt}{$\scriptstyle #2$}}$}}

\newcommand{\nen}{\mbox{$\not \in$}}

\newcommand{\Spec}{\mbox{$\mathrm {Spec\,}$}} 
\newcommand{\sub}{\mbox{$\subseteq$}}

\newcommand{\bcap}{\mbox{$\bigcap$}}
\newcommand{\bcup}{\mbox{$\bigcup$}}
\newcommand{\sth}{\, \vert \,}
\newcommand{\0}{\mbox{$\emptyset$}} 

\newcommand{\Ga}{\mbox{$\Gamma$}}
\newcommand{\ga}{\mbox{$\gamma$}}

\newcommand{\ka}{\mbox{$\kappa$}}

\newcommand{\si}{\mbox{$\sigma$}}

\newcommand{\fh}{\raisebox{.4ex}{\mbox{$\varphi$}}}

\newcommand{\De}{\mbox{$\Delta$}}

\newcommand{\al}{\mbox{$\alpha$}}
\newcommand{\bt}{\mbox{$\beta$}}

\newcommand{\cB}{\mbox{$\cal B$}}
\newcommand{\cC}{\mbox{$\cal C$}}

\newcommand{\cF}{\mbox{$\cal F$}}

\newcommand{\cD}{\mbox{$\cal D$}}
\newcommand{\cH}{\mbox{$\cal H$}}

\newcommand{\cG}{\mbox{$\cal G$}}

\newcommand{\cl}[1]{\mbox{$\overline {#1}$}}

\newcommand{\spez}{\,\mbox{$\rightsquigarrow$\,}}
\newcommand{\spezz}[1]{{\underaccent{#1}{\spez}}}%
\newcommand{\ex}{\mbox{$\exists$}}
\newcommand{\fa}{\mbox{$\forall$}}
\newcommand{\w}{\mbox{$\wedge$}} 
\newcommand{\bo}{\mbox{$\bot$}} 

\newcommand{\Ra}{\mbox{$\Rightarrow$}}
\newcommand{\Lra}{\mbox{$\Leftrightarrow$}}

\newcommand{\La}{\mbox{$\Leftarrow$}}

\newcommand{\und}[1]{\raisebox{-.2ex}{\underline{\raisebox{.2ex}{#1}}}}  
\newcommand{\y}[1]{\mbox{$#1$}}

\newcommand{\bc}{\begin{center}}
\newcommand{\ec}{\end{center}}

\newcommand{\lbr}{\linebreak}

\newcommand{\fm}[1]{\mbox{$\langle \, #1 \,  \rangle$}}

\newcommand{\td}[1]{\mbox{$\widetilde{#1}$}}
\newcommand{\h}[1]{\mbox{$\widehat{#1}$}}

\newtheorem{Th}{Theorem}[section]
\newtheorem{Co}[Th]{Corollary}
\newtheorem{Df}[Th]{Definition}
\newtheorem{Pro}[Th]{Proposition}
\newtheorem{Le}[Th]{Lemma}
\newtheorem{Exa}[Th]{Example}
\newtheorem{Exas}[Th]{Examples}
\newtheorem{Rem}[Th]{Remark}
\newtheorem{Dfre}[Th]{Definition and Remarks}
\newtheorem{Dfno}[Th]{Definition and Notation}
\newtheorem{Rems}[Th]{Remarks}
\newtheorem{Remno}[Th]{Remarks and Notation}
\newtheorem{Fa}[Th]{Fact}
\newtheorem{No}[Th]{Notation}
\newtheorem{Prel}[Th]{Preliminaries}
\newtheorem{Prelnot}[Th]{Preliminaries and Notation}
\newtheorem{Ct}[Th]{\hspace*{-4pt}}

\newcommand{\bre}{\begin{Rem} \em}
\newcommand{\bdfr}{\begin{Dfre}}
\newcommand{\edfr}{\end{Dfre}}
\newcommand{\bdfn}{\begin{Dfno} \em}
\newcommand{\edfn}{\end{Dfno}}
\newcommand{\bdf}{\begin{Df} \em}
\newcommand{\edf}{\end{Df}}
\newcommand{\bremn}{\begin{Remno} \em}
\newcommand{\eremn}{\end{Remno}}
\newcommand{\bth}{\begin{Th}}
\newcommand{\eth}{\end{Th}}
\newcommand{\bco}{\begin{Co}}
\newcommand{\eco}{\end{Co}}
\newcommand{\ble}{\begin{Le}}
\newcommand{\ele}{\end{Le}}
\newcommand{\bpr}{\begin{Pro}}
\newcommand{\epr}{\end{Pro}}
\newcommand{\bex}{\begin{Exa} \em}
\newcommand{\eex}{\end{Exa}}
\newcommand{\bexs}{\begin{Exas} \em}
\newcommand{\eexs}{\end{Exas}}
\newcommand{\ere}{\end{Rem}}
\newcommand{\bres}{\begin{Rems} \em}
\newcommand{\eres}{\end{Rems}}
\newcommand{\bfa}{\begin{Fa}}
\newcommand{\efa}{\end{Fa}}
\newcommand{\bno}{\begin{No} \em}
\newcommand{\eno}{\end{No}}
\newcommand{\bprel}{\begin{Prel} \em}
\newcommand{\eprel}{\end{Prel}}
\newcommand{\bprelnot}{\begin{Prelnot} \em}
\newcommand{\eprelnot}{\end{Prelnot}}
\newcommand{\bct}{\begin{Ct} \em}
\newcommand{\ect}{\end{Ct}}

\newcommand{\llab}[1]{\ar@{}[l]|-<<{\txt{\footnotesize #1}}} 
\newcommand{\sllab}[1]{\ar@{}[l]|-<<{\txt{\scriptsize #1}}}  
\newcommand{\ulab}[1]{\ar@{}[u]|-<<{\txt{\footnotesize #1}}} 
\newcommand{\dlab}[1]{\ar@{}[d]|-<<{\txt{\footnotesize #1}}} 
\newcommand{\rlab}[1]{\ar@{}[r]|-<<<{\txt{\footnotesize #1}}} 
\newcommand{\srlab}[1]{\ar@{}[r]|-<<{\txt{\scriptsize #1}}}  
\newcommand{\lulab}[1]{\ar@{}[l]_<<{\txt{\footnotesize #1}}} 
\newcommand{\rulab}[1]{\ar@{}[r]^<<{\txt{\footnotesize #1}}} 
\newcommand{\ldlab}[1]{\ar@{}[l]^<<<{\txt{\footnotesize #1}}} 
\newcommand{\rdlab}[1]{\ar@{}[r]_>>>>{\txt{\footnotesize #1}}} 
\newcommand{\edge}[1]{\ar@{-}[#1]} 

\def\NewFont#1#2#3#4#5{%
\expandafter\font\csname #1display\endcsname =#1 at #2%
\expandafter\font\csname #1normal\endcsname =#1 at #3%
\expandafter\font\csname #1script\endcsname =#1 at #4%
\expandafter\font\csname #1scriptscript\endcsname =#1 at #5%
}
\def\NewFontLetter#1#2{{%
\mathchoice%
{{\expandafter\hbox{\csname #1display\endcsname\char"#2}}}%
{{\expandafter\hbox{\csname #1normal\endcsname\char"#2}}}%
{{\expandafter\hbox{\csname #1script\endcsname\char"#2}}}%
{{\expandafter\hbox{\csname #1scriptscript\endcsname\char"#2}}}%
}}%
\NewFont{matha10}{9.00truept}{8.00truept}{7.00pt}{6.00pt}

\ifx\undefined\fontshape   
\else
\fi

\begin{document}

\title{Fans in the Theory of Real Semigroups \\
II. Combinatorial Theory}

\vspace{-0.3cm}
\author{M. Dickmann \and A. Petrovich}

\date{March 2017} 

\maketitle

\pagestyle{myheadings}

\vspace*{-1cm}

\begin{abstract}
In \cite{DP5a} we introduced the notion of fan in the categories of real semigoups and their dual abstract real spectra 
and developed the algebraic theory of these structures. In this paper we develop the combinatorial theory of ARS-fans, 
i.e., fans in the dual category of abstract real spectra. Every ARS is a spectral space and hence carries a natural partial 
order called the {\it specialization partial order}. Our main result shows that the isomorphism type of a finite fan in the 
category ARS is entirely determined by its order of specialization. The main tools used to prove this result are: (1) Crucial 
use of the theory of {\it ternary semigroups}, a class of semigroups underlying that of RSs; (2) Every ARS-fan is a disjoint 
union of abstract order spaces (called {\it levels}); (3) Every level carries  a natural involution of abstract order spaces, 
and (4) The notion of a {\it standard generating system}, a combinatorial tool replacing, in the context of ARSs, the 
(absent) tools of combinatorial geometry (matroid theory) employed in the cases of fields and of abstract order spaces.
\end{abstract}

\vspace*{-0.9cm} 

\section*{Introduction} \label{fan-intro} 

\vspace*{-0.4cm}

In \cite{DP5a} we introduced a notion of {\it fan} in each of the dual categories {\bf{RS}} and {\bf{ARS}}
of real semigroups, and of abstract real spectra, dubbed, respectively, {\it RS-fans} and {\it ARS-fans}. 
The emphasis in \cite{DP5a} was on the algebraic theory of RS-fans. The present paper, 
a continuation of \cite{DP5a}, is devoted to develop the combinatorial theory of ARS-fans, i.e., fans
in the category of abstract real spectra.
\vspace*{-0.15cm}

The Introduction to \cite{DP5a} gives an account of the role of fans in the theories of preordered fields,
of quadratic forms, and in real algebraic and analytic geometry.
\vspace*{-0.15cm}

Our main result in this paper is Theorem \ref{isothffans}, showing that the isomorphism type of a finite 
ARS-fan (in the category {\bf{ARS}}) is entirely determined by its order of specialization as a spectral 
space\,\footnote{\,For a general reference on spectral spaces, see \cite{DST}.}. 
The proof of this result relies on a combinatorial machinery that we set up in \S\S\,\ref{levelsARS},\ref{involut}. 
This machinery also gives detailed information on the structure of ARS-fans under their order of specialization.

\vspace*{-0.15cm}

In Section \ref{levelsARS} we introduce the notion of {\it level} of an ARS-fan $(X,F)$. Levels are the pieces 
$\se L I = \{h \in X \sth h^{-1}[0] =  I\}$ of the partition\vspace{0.05cm} of the character space $X$ of the
real semigroup $F$, determined by the ideals $I$ of $F$. Since $F$ is a RS-fan, its ideals are necessarily 
prime and saturated (\cite{DP5a}, Prop. 1.6\,(4), Cor. 3.10\,(1)), and the family of them is totally ordered 
under inclusion (\cite{DP5a}, Fact 1.4). By Proposition 5.11 of \cite{DP5a}, each level is an abstract space 
of orders and therefore (by results from \cite{D1}, \cite{D2} and \cite{Li}) possesses a structure of 
combinatorial geometry (matroid). Further, there exist canonical AOS-morphisms linking each level to 
any level determined by a larger ideal (i.e., a ``higher'' level); Proposition \ref{RSG+levelMorph}\,(2). 
\vspace*{-0.15cm}

In the next \S\,\ref{involut} we exploit the combinatorial geometric structure of the levels to investigate the
fine structure of ARS-fans. In Theorem \ref{propsinvol} we show that multiplication of a character
of level $\se L I$ by any pair of elements $\se g 1, \se g 2 \in X$ so that $Z(\se g i) := {\se g i}^{-1}[0]\,
\sub\, I\; (i = 1,2)$ defines an involution of the AOS $\se L I$. These involutions are compatible
with the order of specialization between levels induced by inclusion of the determining ideals  
(\ref{propsinvol}\,(e)). Further, we prove that these involutions permute certain AOS-subfans of the levels 
defined by combinatorial conditions (Propositions \ref{cardsetsS+C} and \ref{permsetsS+C}). Altogether, 
the results proved in this section show that the order structure of ARS-fans is subject to strong
constraints, illustrated in \ref{impossibles}.
\vspace*{-0.15cm}

To prove the isomorphism Theorem \ref{isothffans}, the combinatorial machinery mentioned above
is used together with the notion of a {\it standard generating system} introduced in \ref{stgen}.
This notion is a substitute for the combinatorial geometric notions existing in the context of AOSs, but
absent in that of ARSs.

\noi {\bf Preliminaries.} \label{prelims} For easy reference we state, without proof, the following simple 
facts proved in \cite{DP5a} and frequently used below. The axioms defining the notion of a {\it ternary
semigroup} (abbreviated TS) appear in \cite{DP5a}, Def. 1.1, and \cite{DP1}, \S\,1, p. 100; $\se X T$
denotes the set of TS-homomorphisms of a TS, $T$, into the TS ${\bf 3} = \{1,0,-1\}$ (the TS-{\it characters}
of $T$).
\vspace*{-0.1cm}

The first Lemma gives several characterizations of the specialization order of the spectral topology on the
character set of a ternary semigroup.
\vspace{-0.45cm}

\ble \label{char-specializ} Let $T$ be a TS, and let $g,h \in \se X T$. The following are equivalent:
\vspace{-0.15cm}

\nhp $(1)$ $g \spez \, h$ $($i.e., \y h is an specialization of \y g$)$.
\vspace{-0.15cm}

\nhp $(2)$ $h^{-1}[1]\, \sub \, g^{-1}[1]$\, $($equivalently,\, $h^{-1}[-1]\, \sub \, g^{-1}[-1])$.
\vspace{-0.15cm}

\nhp $(3)$ $g^{-1}[\{0,1\}]\, \sub \, h^{-1}[\{0,1\}]$.
\vspace{-0.15cm}

\nhp $(4)$ $Z(g)\, \sub \, Z(h)$ and $\fa\, a \in G\, (a \not\in Z(h) \;\; \Ra \;\; g(a) = h(a))$.
\vspace{-0.15cm}

\nhp $(5)$ $h = h^2g$ $($equivalently, $h^2 = hg)$. \hfl $\Box$
\ele 

\vspace{-0.45cm}

We also register the following algebraic characterizations of inclusion 
and equality of zero-sets of elements of $\se X T$.

\vspace{-0.45cm}

\ble \label{char-zeroset} Let $T$ be a TS, and let $u,g,h \in \se X T$. Then:
\vspace{-0.15cm}

\nhp $(1)$ $Z(g)\, \sub\, Z(h)\; \Lra\; h = h g^2$.
\vspace{-0.15cm}

\nhp $(2)$ $Z(g) = Z(h)\; \Lra\; g^2 = h^2$.
\vspace{-0.15cm}

\nhp $(3)$ If\, $u \spez \, g, h$, then\, $Z(g)\, \sub \, Z(h)$\, if and only if\, 
$g \spez \, h$. \hfl $\Box$
\ele

\vspace{-0.75cm}

\bpr \label{charfan:(1)implies(2.ii)} Let $G$ be a RS-fan. Then: 

\vspace{-0.15cm}
\nhp $(1)$ For all elements $g,\,h \in \se X G$ such that $g \spez\, h$ $($hence
$Z(g) \, \sub \, Z(h))$ and every ideal \y I such that $Z(g) \, \sub \, I \,
\sub \, Z(h)$ there is $f \in \se X G$ such that $g \spez f \spez\, h$ and $Z(f) = I$.
\vspace{-0.1cm}

\nhp $(2)$ For every \, $g \in \se X G$ and every ideal \,$I \supseteq Z(g)$ there is 
a $($necessarily unique$)$ $f \in \se X G$  such that $g \spez f$ and
$Z(f) = I$. 
\vspace{-0.1cm}

\nhp $(3)$ For every ideal \y I of \y F there is an $f \in \se X G$ such that
$Z(f) = I$. \hfl $\Box$
\epr
\vspace{-0.5cm}

\section{Levels of a ARS-fan} \label{levelsARS}

\vspace{-0.3cm}

The saturated prime ideals of a real semigroup induce a partition of its
character space. The pieces are called {\it levels}\,: the level corresponding to
a saturated prime ideal $I$ of $G$ is the set of all $g \in \se X G$ such that
$Z(g) = I$. In the case of RS-fans, (proper) ideals ---automatically prime 
 and saturated (\cite{DP5a}, Prop. 1.6\,(4), Cor. 3.10\,(1))--- are totally ordered 
 under inclusion (\cite{DP5a}, Fact 1.4), a fact of much help in studying the relationship 
 between their levels. This notion, together with that of a {\it connected component} 
 (\ref{components}), will be the main technical tools employed in the analysis of the 
 fine structure of ARS-fans carried out in this paper.
 \vspace{-0.15cm}

Proposition 5.11 in \cite{DP5a} shows that the levels of an ARS-fan have a canonical structure of 
AOS-fans (\ref{levels}\,(a)), that is, of fans in the category of abstract order spaces 
(cf. \cite{M}, \S\,3.1, pp. 37 ff.). We prove (Proposition \ref{RSG+levelMorph}) that inclusion 
of ideals induces AOS-morphisms between the corresponding levels (cf. \ref{levels}\,(c)). As an
application we prove (Corollary \ref{cardffans}) that the cardinality of a finite RS-fan, $F$, and 
that of its character space, $\se X F$, are related by the identity ${\mathrm {card}}\,(F) = 
2 \cdot {\mathrm {card}}\,(\se X F) + 1$, an analog for RSs of a result known to hold for 
reduced special groups.  \hfl $\Box$
\vspace{-0.3cm}

\bprelnot \label{levels} (a) Given a real semigroup \y G and a saturated prime ideal 
\y I of \y G, we denote by \se G I the RSG $(G/I) \setminus \{ \pi(0) \}$, cf. \cite{DP5a}, Prop. 5.11. 
The congruence of $G$ determined by the set of characters $\se {\cH} I = \{h \in \se X G \sth Z(h) = I\}$ is 
denoted by $\se {\sim} I$, cf. \cite{DP5a}, \S\,5.C. Every character $h \in \se {\cH} I$ induces a map\, 
$\h h : \se G I \, \lra \, \{\, \pm 1\}$ defined by\, $\h h \, \com \, \se {\pi} I = h$. The correspondence 
$h\, \mapsto \h h$ is a bijection between the set\, $\se L I(G) = \se {\cH} I$\, and the space of 
orders \se X {G_I}\, of\, \se G I\,. ($\se L I$ stands for ``$I$-th level"; see item (b.ii) below). Thus, we can identify 
the set $\se L I(G) \; \sub \, \se X G$ with the AOS $(\se X {G_I}, \se G I)$. We shall systematically use 
this identification in the sequel, and unambiguously refer to the \und{AOS structure} of the set 
\se L I(\y G). In case \y G is a RS-fan, \cite{DP5a}, Prop. 5.11, shows that \se L I(\y G) is an \und{AOS-fan}.

\noi (b) \label{depthlev} Let \y F be a RS-fan.
\vspace{-0.15cm}

\nhp (i)\;\; We denote by\: $\Spec (F)$\, the set of all proper ideals of \y F.
\vspace{-0.15cm}

\nhp (ii)\; For $I \in \Spec (F)$ the set $\se L I(F) = \{ h \in \se X F
\, \vert \, Z(h) = I \}$ is called the {\bf \y I-th level} of $\se X F$\,.
\vspace{-0.15cm}

\nhp (iii) For $f \in \, \se X F$, the {\bf depth} of \y f, denoted $d(f)$, 
is the order type of the set $\{ g \in \se X F\, \vert \, f \spez \, g \}$ 
under the order of specialization.\index{sub}{depth}\index{sub}{fan!depth of
element} (Since $(\se X F, \spez)$ is a root-system, the order $\spez$ is
total on this set.)
\vspace{-0.15cm}

\nhp (iv)\, For $I \in \Spec (F)$, the order type of the set 
$\{ J \in \Spec (F)\, \vert \, J \supseteq \, I \}$ under the 
$($total\,$)$ order of inclusion will be called the {\bf depth} of \y I, 
denoted $d(I)$.
\vspace{-0.15cm}

\nhp (v)\; The {\bf length} of \se X F, denoted $\ell(\se X F)$, is the order 
type of the $($totally ordered\,$)$ set\, $\Spec (F)$.

\vspace{-0.1cm}

\noi (c) (AOS- and ARS-morphisms; \cite{M}, \S\,2, pp. 23-24, and \S\,6, p. 103)
\vspace{-0.15cm}

\noi (i) Let $(X,G), (Y,H)$ be ARS's. A map $F: X\, \lra\, Y$ is an {\bf ARS-morphism} iff for 
all\, $a \in H$\, there is\, $b \in G$\, so that\, $\h a\; \com\, F = \h b$. Here, for $x \in G$,
$\h x: X\, \lra\, {\bf 3}$ denotes the map ``evaluation at $x$'': $\h x(\si) := \si(x)$, for 
$\si \in X$, and similarly for $H$.
\vspace{-0.15cm}

\noi (ii) The definition of an {\bf AOS-morphism} is similar, with $(X,G), (Y,H)$ AOS's, and the
evaluation maps taking values in $\{\pm 1\}$.
\vspace{-0.15cm}

\noi (d) If $f:G\, \lra\, H$ is a RS-morphism\vspace{0.08cm} (resp. RSG-morphism), the dual map 
$f^*: \se X H\, \lra\, \se X G$ defined by $f^*(\ga) := \ga\, \com\, f$\, for $\ga \in \se X H$, 
is an ARS-morphism (resp., AOS-morphism). \hfl $\Box$ 

\vspace{-0.1cm}

\nhp {\bf Remarks.} (a) Clearly, the union and the intersection of an 
inclusion chain of (proper) prime ideals in any ternary semigroup is a 
(proper) prime ideal. In particular, if \y F is a fan, the totally ordered 
set $(\Spec(F), \sub)$ is (Dedekind) {\bf complete}.  \hfl $\Box$
\eprelnot

\vspace{-0.8cm}

\bpr \label{RSG+levelMorph} Let \y F be a RS-fan and let $I \, \sub \, J$ be ideals 
of \y F. With notation as in \ref{levels}, 
\vspace{-0.15cm}

\noi $(1)$ The rule $a/J \, \mapsto \, a/I$\; $(a \in F \setminus J)$ defines a 
homomorphism of special groups $\se {\iota} {JI}: \se F J \, \lra \, \se F I$.

\vspace{-0.15cm}

\noi $(2)$ The map $\se {\ka} {IJ}: \se L I(F)\, \lra \, \se L J(F)$ assigning 
to\vspace{0.05cm} each $g \in \se L I(F)$ the unique element $h \in \se L J(F)$ 
such that $g \spez \, h$ is an AOS-morphism.
\epr

\vspace{-0.4cm}

\nhp {\bf Proof.} (1) (a) $\se {\iota} {JI}$ is well-defined.
\vspace{-0.15cm}

\nhp We must show: $a,b \, \in F \setminus J  \;\w \; a \, \se {\sim} J \, b
\;\,\Ra \;\, a \, \se {\sim} I\, b$. Since $I \, \sub \, J$, this is clear from
Lemma 5.10 of \cite{DP5a} which states that, for an ideal $K$ of $F$
and $a,b \in F \setminus K,\;\, a \, \se {\sim} J \, b\; \Lra\; \ex z \not\in K\, (az = bz)$. 
\vspace{-0.15cm}

Clearly, we have:
\vspace{-0.15cm}

\nhp (b) $\se {\iota} {JI}$ is a group homomorphism sending $-1/J$ to $-1/I$.
\vspace{-0.15cm}

\nhp Since \se F J is a RSG-fan, \se {\iota} {JI} is automatically a homomorphism
of special groups. 
\vspace{-0.15cm}

\nhp (2) By (1) and \ref{levels}\,(d), the map $\se {\iota^*} {JI}: \se X {F_I}\, \lra \, \se X {F_J}$ 
dual to $\se {\iota} {JI}$ is an AOS-morphism. The map $\se {\ka} {IJ}$ is \, 
$\se {\ka} {IJ} = ({\se {\fh} {\!J}})^{-1} \circ \se {\iota^*} {JI} \circ \se {\fh} {\!I}$, 
where $\se {\fh} {\!I}$ denotes the bijection $g \mapsto \h g$ ($g \in \se L I(F)$), 
identifying $\se L I(F)$ with \se X {F_I} (\ref{levels}\,(a)), and similarly for 
\se L J(\y F). It only remains to prove $g \spez \, \se {\ka} {IJ}(g)$, 
for $g \in \se L I(F)$. To ease notation, write $h = \se {\ka} {IJ}(g)$.
According to\vspace{0.05cm} Lemma \ref{char-specializ}\,(4) we must show $Z(g) \, \sub \, Z(h)$ 
and $a \not\in Z(h) \; \Ra \; g(a) = h(a)$. The inclusion of zero-sets is $I \, \sub \, J$. 
Let $a \not\in Z(h) = J$. Since:
\vspace{-0.1cm}

\nhp \hfl $\se {\fh} {\!J}(h) = \se {\fh} {\!J} (\se {\ka} {IJ}(g)) = \se {\iota^*}
{JI}(\se {\fh} {\!I}(g)) = \se {\fh} {\!I}(g)\, \com\, \se {\iota} {JI}$\,, \; $\se
{\fh} {\!I}(g)(a/I) = g(a)$\; and\; $\se {\fh} {\!J}(h)(a/J) = h(a)$, \hfl
\vspace{-0.1cm}

\nhp (cf. \ref{levels}), we get,
\vspace{-0.2cm}

\nhp \hfl $h(a)= (\se {\fh} {\!I}(g))(\se {\iota} {JI}(a/J)) = \se {\fh} {\!I}(g)(a/I) = g(a)$, \hfl
\vspace{-0.2cm}

\nhp as required. \hfl $\Box$
\vspace{-0.1cm}

Next we prove that the depth of an ideal in a fan is the same as the depth
of any element in the corresponding level; in particular, elements of the same
depth belong to the same level.

\vspace{-0.4cm}

\bpr \label{depthequ} Let \y F be a RS-fan. For $f \in \se X F$ we have $d(f) =
d(Z(f))$; equivalently, the sets $\{ g \in \, \se X F \, \vert \, f \spez \, g\}$ 
$($ordered under specialization\,$)$ and $\{ J \in$ {\em Spec}$(F)\, \vert \, J \supseteq Z(f) \}$ 
$($ordered under inclusion\,$)$ are order-isomorphic.
\epr

\vspace{-0.4cm}

\nhp {\bf Proof.} To ease notation, set $f\! \uparrow \; = \{\, g \in \, \se X F
\, \vert \, f \spez\, g\}$ and $I\!\uparrow \; = \{ J \in \textrm{Spec}(F)\,
\vert \, J \supseteq I \}$ \lbr
$(I \in \textrm{Spec}(F))$. The required order isomorphism is the map 
$Z: f\!\uparrow \, \lra \; Z(f)\!\uparrow$\; assigning to each\, 
$g \in f\!\uparrow$\, its zero-set. That
\vspace{-0.15cm}

\nhp (a)\; \y Z is increasing, \;\; and \;\; (b)\; \y Z is surjective,
\vspace{-0.15cm}

\nhp is clear, from $g \spez \, h \;\, \Ra \;\, Z(g) \, \sub \, Z(h)$ and 
Proposition \ref{charfan:(1)implies(2.ii)}\,(2), respectively. That 
\vspace{-0.15cm}

\nhp (c) \y Z is injective.
\vspace{-0.15cm}

\noi follows from \ref{char-zeroset}\,(3).   \hfl $\Box$

A trivial variant of the proof of \ref{depthequ} gives:

\vspace{-0.4cm}

\bpr \label{isointerv} Let \y F be a RS-fan. Given \,$f_1, f_2 \in \se X F$ such
that $f_1 \spez\, f_2$, the intervals $\{\, g \in \se X F \, \vert \, f_1 \spez
\, g \spez \, f_2 \}$ $($ordered under specialization\,$)$ and $\{ J \in
Spec(F)\, \vert \, Z(f_1) \, \sub \, J \, \sub \, Z(f_2) \}$ $($ordered under
inclusion\,$)$ are order-isomorphic. \hfl $\Box$
\epr

\vspace{-0.3cm}

The results in the next two Lemmas will be frequently used in this and in subsequent sections.
\vspace{-0.4cm}

\ble \label{prods} Let $G$ be a RS and let $\se g 1, \dots ,\, \se g r, h \in \se X G$ be so that\, 
$\bcup_{i=1}^r Z(\se g i) \, \sub \, Z(h)$. For $i = 1, \dots, r$, let $\se f i \in \se X G$ 
be such that $\se g i \spez \, \se f i$ and \,$Z(\se g i) \, \sub \, Z(\se f i) \, \sub \, Z(h)$. 
Then,
\vspace{-0.15cm}

\nhp $(*)$ \;\;\; $h \cdot\se g 1 \cdot \dots\, \cdot \se g r = h \cdot \se f 1 \cdot \dots\, 
\cdot \se f r$.
\ele

\vspace{-0.45cm}

\nhp \und{Note}. The products in (*) may not be in $\se X G$.

\nhp {\bf Proof.} Obviously, (*) holds whenever $x \in Z(h)$. If $x \not\in Z(h)$, from 
the assumptions we get $x \not\in \, \bcup_{i=1}^r Z(\se g i)$\, and\, 
$x \not\in \;\bcup_{i=1}^r Z(\se f i)$. Since $\se g i \spez \se f i$, we get 
$\se g i(x) = \se f i(x)$ for $i = 1, \dots,\, r$ (Lemma \ref{char-specializ}\,(4)), 
and (*) follows. \hfl $\Box$

\vspace{-0.4cm}

\ble \label{moreprods} Let $F$ be a RS-fan. Then,
\vspace{-0.15cm}

\noi $(a)$\, For\, $i = 1, \dots,\, r$, with \y r \und{odd}, let \,
$\se g i, \se h i \in \se X {\!F}$ be such that \,$\se g i \spez \, \se h i$. Then,\, 
$\se g 1\cdot \dots\, \cdot\se g r \, \spez$ \, $\se h 1\cdot \dots \, \cdot\se h r$.
\vspace{-0.15cm}

\nhp $(b)$\, Let $\se h 1, \se h 2, f, g, k \in \se X {\!F}$ be such that $f,\,g\, \spez \, \se h 1, \;
k\, \spez \, \se h 2$, and $Z(\se h 1) \, \sub\, Z(\se h 2)$. Then,\lbr 
$fg\,k\, \spez \, \se h 2$.
\ele

\vspace{-0.45cm}

\nhp \und{Note}. Here the products are in \se X {\!F} as the number of factors is odd.

\nhp {\bf Proof.} (a) For\, $i = 1, \dots,\, r$ we have $\se {h^2} i = \se h i \se g i$ 
(Lemma \ref{char-specializ}\,(5)). Multiplying these equalities termwise gives 
$(\se h 1 \cdot \dots \, \cdot \se h r)^2 = (\se h 1 \cdot \dots \, \cdot \se h r)
(\se g 1 \cdot \dots \, \cdot \se g r)$, which proves the assertion.
\vspace{-0.15cm}

\nhp (b) By Lemma \ref{char-specializ} we must prove $\se {h^2} 2 = \se h 2 (f g\, k)$. Obviously, 
this\vspace{0.05cm} equality holds at every $x \in Z(\se h 2)$. If $x \not\in Z(\se h 2)$, then 
$x \not\in Z(\se h 1)$, and 
$f,\,g \spez \, \se h 1$ implies $\se h 1(x) = f(x) = g(x) \neq 0$; also $k \spez \, \se h 2$ implies 
$\se h 2(x) = k(x) \neq 0$, whence $f(x)g(x) = 1$ and $\se h 2(x)k(x) = 1$. This yields\, 
$(\se h 2 f g\, k)(x) = (f(x)g(x))(\se h 2(x)k(x)) = 1$. On the other hand, $(\se h 2(x))^2 =1$, 
proving that the required identity holds at $x \not\in Z(\se h 2)$\, as well.  \hfl $\Box$

Our last result in this section, Corollary \ref{cardffans}, shows that if $F$ is a \und{finite} RS-fan
and $\se X F$ its character space, then ${\mathrm {card}}\,(F) = 2 \cdot {\mathrm {card}}\,(\se X F) + 1$. 
This identity is the analog of a well-known result relating the cardinalities of a finite RSG-fan 
and its space of orders (\cite{ABR}, p. 75). The result follows from a more general observation, 
valid for RS-fans of arbitrary cardinality.
\vspace{-0.35cm}

\bpr \label{IsoQuotideals} Let $I \subset J$ be consecutive ideals of a RS-fan $F$ $($with, possibly, $J = F)$. Then,

\vspace{-0.15cm}

\noi $($i\,$)$\;\, Under product induced by $F$, $J \setminus I$ is a group of exponent $2$ with unit $x^2$
for any $x \in J \setminus I$ $($and distinguished element $-1 = -x^2)$.
\vspace{-0.15cm}

\noi $($ii\,$)$ The restriction of the quotient map $\se {\pi} I\, \lceil\,(J \setminus I): J \setminus I\, 
\lra\, \se F I = F/I \setminus \{\se {\pi} I(0)\}$ is a group isomorphism preseving the distinguished
element $-1$.
\epr
\vspace{-0.38cm}

\noi {\bf Proof.} (i) Since $I$ is prime, $J \setminus I$ is closed under product. Given 
$x, y \in J \setminus I$, we must prove $x^2 = y^2$ (which implies $x^2y = y^3 = y$). By 
the separation theorem for TSs (\cite{DP1}, Thm. 1.9, pp. 103-104) it suffices to show that 
$h(x^2) = h(y^2)$ for all $h \in \se X F$. If $J\, \sub\, Z(h)$, then $h(x^2) = h(y^2) = 0$. 
If $Z(h)\, \sub\, I$, then $h(x), h(y) \neq 0$, whence $h(x^2) = h(y^2) = 1$.
\vspace{-0.15cm}

\noi (ii) Clearly, $\se {\pi} I(x) \neq \se {\pi} I(0)$, i.e., $\se {\pi} I(x) \in \se F I$,
for all $x \in J \setminus I$, and $\se {\pi} I$ preserves product.
\vspace{-0.15cm}

\noi --- $\se {\pi} I\, \lceil\,(J \setminus I)$ is injective.
\vspace{-0.15cm}

\noi Suppose $\se {\pi} I(x) = \se {\pi} I(y)$, i.e., $x\: \se {\sim} I\, y$, with 
$x, y \in J \setminus I$. By\vspace{0.05cm} \cite{DP5a}, Lemma 5.10 (cf. proof of 
\ref{RSG+levelMorph}), $xz = yz$ for some $z \not\in I$. 
To prove $x = y$, let $h \in \se X F$. If $J\, \sub\, Z(h)$, then $h(x) = h(y) = 0$. If 
$Z(h)\, \sub\, I$, then $h(z) \neq 0$, and we get $h(x) = h(y)$.
\vspace{-0.15cm}

\noi --- $\se {\pi} I(x^2) = \se {\pi} I(1)$, for $x \in J \setminus I$.
\vspace{-0.15cm}

\noi Clear, for $Z(h) = I$ implies $h(x^2) = 1$. In particular, $\se {\pi} I$ preserves $-1$. 
\vspace{-0.15cm}

\noi --- $\se {\pi} I\, \lceil\,(J \setminus I)$ is onto $\se F I$.
\vspace{-0.15cm}

\noi Let $p \in \se F I$\,; then, $p = \se {\pi} I(q)$ with $q \not\in I$. Taking
$z \in J \setminus I$, we have $qz^2 \in J \setminus I$, whence 
$\se {\pi} I(qz^2) = \se {\pi} I(q) \se {\pi} I(z^2) = \se {\pi} I(q) \se {\pi} I(1) =
\se {\pi} I(q) = p$. \hfl $\Box$

\vspace{-0.4cm}

\bno \label{finfannotat} Let \y F be a finite RS-fan, and let
\vspace{-0.15cm}

\nhp \hfl $\{ 0 \} = \se I n \subset \, \se I {n-1} \subset \cdots \subset \se I 2 
\subset \, \se I 1  \subset \, F = \se I 0$ \hfl
\vspace{-0.15cm}

\nhp be the set of all its ideals; thus, for $1 \leq d \leq n$, \se I d is the ideal of depth \y d. 
We set $\se F d = \se F {I_d} = (F/\se I d) \setminus \{ \se {\pi} d(0)\}$, where 
$\se {\pi} d: F\, \lra \; F/\se I d$ denotes the canonical quotient map. We 
also write \se L d for $\se L {I_d}$; cf. \ref{depthlev}\,(b).  
\hfl $\Box$
\eno  

\vspace{-0.5cm}

\nhp Clearly, $F \setminus \{ 0 \} = \bcup_{d=1}^n (\se I {d-1} \setminus \se I
d)$ (disjoint union),\vspace{0.08cm} whence, by \ref{IsoQuotideals} we have 
$\mbox{$\mathrm {card}$}\,(F) = \sum_{d=1}^n \mbox{$\mathrm {card}$}\,(\se I {d-1} \setminus \se I d) + 1 = 
\sum_{d=1}^n \mbox{$\mathrm {card}$}\,(\se F d) + 1$. Since the levels partition \se X F, \ref{levels} yields:

\vspace{-0.4cm}

\bpr \label{card1} For any finite RS-fan \y F, \mbox{$\mathrm {card}$}\,$(\se X F) = \sum_{d=1}^n
\mbox{$\mathrm {card}$}\, (\se L d) = \sum_{d=1}^n {\mathrm {card}}\,(\se X {\!F_d})$.\hfl $\Box$
\epr
\vspace{-0.7cm}

\bco \label{cardffans} For a finite RS-fan, \y F, \; ${\mathrm {card}}\,(F) = 
2 \cdot {\mathrm {card}}\,(\se X F) + 1$.
\eco
\vspace{-0.45cm}

\nhp {\bf Proof.} Since the \se F d are finite RSG-fans (\cite{DP5a}, Prop. 5.11), 
we know\vspace{0.05cm} that ${\mathrm {card}}\,(\se F d) =$\lbr
$2 \cdot {\mathrm {card}}\,(\se X {\!F_d})$ for $1 \leq \, d \, \leq \, n$ (see \cite{ABR},
p. 75). The result follows, then, from Proposition \ref{IsoQuotideals} and the preceding 
cardinality identities. \hfl $\Box$

\vspace{-0.4cm}

\section{Involutions of ARS-fans}\label{involut}

\vspace{-0.2cm}
\bno \label{setsS+C} In addition to the notation introduced in Definition 
\ref{depthlev}, for $J \, \sub \, I$ in Spec(\y F) we define the sets:
\vspace{-0.15cm}

\nhp \hspace{0.5cm} $\se {S^I} {\!J} = \{\, h \in \se L I \, \vert \, \ex g \in \se X {\!F}\,
(\,g \spez\, h \; \w \; Z(g) = J\,) \}.$
\vspace{-0.2cm}

\nhp \hspace{0.5cm} $\se {C^I} {\!J} = \{\, h \in \se {S^I} {\!J} \, \vert \, 
\fa g' \in \se X {\!F}\,(\,g' \spez\, h \; \Ra \; J \, \sub \, Z(g'))\}.$
\vspace{-0.15cm}

\nhp That is,\, $\se {S^I} {\!J}$ consists of those elements of level \y I having 
predecessors of level \y J \und{or lower} in the specialization partial order; 
\,$\se {C^I} {\!J}$ is the set of elements in \se L I having predecessors at level 
\y J \und{but not lower}. \hfl $\Box$
\eno

\vspace{-0.7cm}

\bres \label{propsetsS+C} (i) For $I \in$ Spec(\y F), $\se {S^I} {\{ 0 \}} = \se {C^I} {\{ 0 \}}$, 
and $\se {S^I} I$ = \se L I. (Recall that $\{ 0 \}$ is the smallest member of Spec(\y F), i.e., 
the zero-set of the lowest level of $\se X {\!F}$.)
\vspace{-0.15cm}

\nhp (ii) \;\,For $J \, \sub \, I$ in Spec(\y F),\, $\se {S^I} {\!J} \neq \0$.
\vspace{-0.2cm}

\nhp \und{Proof}. Let $g \in \se X {\!F}$ be such that $Z(g) = J$ (exists by Proposition 
\ref{charfan:(1)implies(2.ii)}\,(3)). If \y h is the unique $\spez$\,-successor of \y g 
of level \y I (Proposition \ref{charfan:(1)implies(2.ii)}\,(2)), then $h \in \se {S^I} {\!J}$.
\vspace{-0.2cm}

\nhp (iii)\, For $J \, \sub \, I$ in Spec(\y F),\, $\se {S^I} {\!J} = {\mathrm {Im}}\,(\se{\ka}{JI})$, where
$\se {\ka} {JI}: \se L {\!J}(F)\, \lra \, \se L I(F)$ is the AOS-morphism defined in Proposition 
\ref{RSG+levelMorph}\,(2).
\vspace{-0.15cm}

\nhp (iv) \,For $J \, \sub \, I$ in Spec(\y F),\, $\se {S^I} {\!J} \supseteq\ \bcup \:
\{\:\se {C^I} {\!J'} \, \vert \, J' \in \, \mbox{Spec}(\y F) \;\, \mbox{and}\;\, J' \, \sub \, J\: \}$. 
(Note that $\se {C^I} {\!J'}$ may be empty for some $J' \, \sub \, J$.)
\vspace{-0.15cm}

\nhp (v) \, For $J \, \sub \, I$ in Spec(\y F),\, $\se {C^I} {\!J} = \se {S^I} {\!J} \setminus 
\bcup \:\{\,\se {S^I} {\!J'} \, \vert \, J' \in \, \mbox{Spec}(\y F) \;\, \mbox{and}\;\, J' \, 
\subset \, J\: \}$.
\vspace{-0.15cm}

\nhp (vi) \,For $J, J' \, \sub \, I$ in Spec(\y F), $J \neq J'$, we have $\se {C^I} {J} \, \cap \,
\se {C^I} {\!J'} = \0$.  \hfl $\Box$
\eres
\vspace{-0.35cm}

In order to render later arguments as transparent as possible, we recall the following simple 
(and well-known) facts about fans in the categories {\bf RSG} and {\bf AOS}.

\vspace{-0.4cm}

\ble \label{fanimage} Let\, $g: H\, \lra\, G$\, be a SG-homomorphism between RSG-fans, and let 
$g^*: (\se X G, G)$ $\lra\, (\se X H, H)$ denote the AOS-morphism dual to \y g $($cf. \ref{levels}\,$(d))$. Then,
\vspace{-0.15cm}

\nhp $(1)$ With representation induced by that of \y H, ${\mathrm {Im}}(g)$ is a RSG-fan, and 
\y G is isomorphic to the extension of\, ${\mathrm {Im}}(g)$ by the exponent-two group 
$\De =  G/\,{\mathrm {Im}}(g)$.
\vspace{-0.15cm}

\nhp $(2)$ $({\mathrm {Im}}(g^*),\, H/\,ker(g))$ is an AOS-fan.
\ele

\vspace{-0.5cm}

\bres \label{ExtRSG+etc} (a) For the definition of extension of a SG by a group of exponent two, 
see \cite{DM1}, Ex. 1.10, p. 10.   
\vspace{-0.15cm}

\nhp (b) By the duality between RSGs and AOSs (\cite{DM1}, Ch. 3), the dual statement holds 
as well: given an AOS-morphism of (AOS-)fans, $\ka: (X,G)\, \lra\, (Y,H)$, the assertions (1) 
and (2) hold with \y g := $\ka^*$ (the SG-morphism dual to \ka), and with $g^* = \ka$. \hfl $\Box$
\eres
\vspace{-0.2cm}

\nhp {\bf Sketch of proof of \ref{fanimage}.} (1) The first assertion is easily checked. For the 
second, ${\mathrm {Im}}(g)$ is a direct summand of the group \y G. Let $pr:\;G\, \lra \;{\mathrm {Im}}(g)$ 
be the projection onto the factor ${\mathrm {Im}}(g)$; $pr$ is a SG-morphism (\y G and ${\mathrm {Im}}(g)$ 
are fans), and is the identity on ${\mathrm {Im}}(g)$. The isomorphism between \y G and 
${\mathrm {Im}}(g)[\De]$ is\; $f(a) = \fm {pr(a),\, a/{\mathrm {Im}}(g)}$, for $a \in\, G$.
\vspace{-0.2cm}

\nhp (2) Recall that $g^*$ is defined by composition,\; $g^*(\si) = \si\, \circ\, g\ (\si \in \se X G)$, see
\ref{levels}\,(d),\, and that\; ${\mathrm {Im}}(g^*)^{\bo} = \bcap\, \{ ker(\ga)\, \vert \, \ga \in\, 
{\mathrm {Im}}(g^*)\} = \bcap\, \{ ker(\si\, \com\: g) \, \vert \, \si \in \se X G \}$. Routine checking from these 
definitions proves that ${\mathrm {Im}}(g^*$) is closed under product of any three members (since 
$\se X G$ is),\vspace{-0.05cm} and that\, ${\mathrm {Im}}(g^*)^{\bo} = ker(g)$\, (since 
$\bcap\, \{ker(\si) \sth \si \in \se X G \} = \{1\}$), whence ${\mathrm {Im}}(g^*)\, \sub\, \se X {H/ker(g)}$.
\vspace{-0.15cm}

Clearly, the map $\cl g : H/ker(g)\, \lra \, {\mathrm {Im}}(g)$ induced by \y g is an SG-isomorphism. 
Thus, we have a commutative diagram of SG-morphisms:

\vspace{-0.45cm}

$$\xymatrix@=3ex{
& H \ar[dr]_-{\pi} \ar[r]^-g & \; {\mathrm {Im}}(g)\; \ar@<1ex>@{^{(}->}[r] & \; G \, \stackrel f {\cong}\, {\mathrm {Im}}(g)[\Delta] \ar@<1ex>[l]^-{pr} \\
& & \; H/ker(g) \ar[u]_{\overline g} & & \\ }
$$

\vspace{-0.2cm}

It only remains to show that ${\mathrm {Im}}(g^*) \, \supseteq \, \se X {H/ker(g)}$. Any 
SG-character $\ga :\, H/ker(g) \, \lra$  $\zt$ can be lifted to a map $\si:\, G \, \lra \, \zt$, 
via the identification of\, \y G with\, ${\mathrm {Im}}(g)[\De]$, as follows: for each $a \in G$ 
there is $b \in H$ such that $pr(a) = g(b)$. We set $\si(a) = \ga(b/ker(g)) = \ga(\pi(b))$. 
In terms of the diagram above, we have: $\si = \ga\, \com\, (\cl g)^{-1}\, \com\, pr$. It follows 
that\, \si\, is a well-defined SG-morphism, i.e., $\si \in \se X G$, and (since $pr\, \com\, g = g$ 
and $(\cl g)^{-1}\, \com\; g = \pi$), $g^*(\si) = \si\, \com\, g = \ga\, \com\, \pi$. \hfl $\Box$

\vspace{-0.15cm}

Lemma \ref{fanimage}, together with \ref{propsetsS+C}\,(iii) and \ref{RSG+levelMorph}\,(2), gives:

\vspace{-0.3cm}

\bco \label{fansS} Let \y F be a RS-fan, and let $J \, \sub \, I$ be in {\em Spec}$(\y F)$. The 
set $\se {S^I} {\!J}$ is an AOS-fan. Indeed, it is a sub-fan of $\se L I(\y F)$, when the latter 
is endowed with its structure of AOS-fan, as indicated in $\ref{levels}$. More generally, if\;
$\cF\, \sub \, \se L {\!J}(\y F)$ is an AOS-fan, the set\, $\se {S^I} {\!J}(\cF) =$ 
$\{\, h \in \se L I \, \vert \,$ $\ex g \in \cF\,(\,g \spez\, h) \}$ is an AOS-subfan of $\se L I(\y F)$. 
\eco

\vspace{-0.3cm}

\nhp {\bf Proof.} The first assertion is a special case of the second (with \cF\, = \se L {\!J}(\y F)). 
For the latter, observe that $\se {S^I} {\!J}(\cF) = \se {\ka} {JI}[\cF] = 
{\mathrm {Im}}(\se {\ka} {JI} \lceil \cF)$ and use Remark \ref{ExtRSG+etc}\,(b). \hfl $\Box$

The following definition will have a crucial role in the sequel:

\vspace{-0.45cm}

\bdf \label{definvol} Let \y F be a RS-fan, let $\se g 1,\se g 2 \in \se X {\!F}$, and fix 
$I \in {\mathrm {Spec}}(F)$\vspace{0.05cm} so that $Z(\se g 1),$  $Z(\se g 2) \, \sub$ \y I. 
We define a map \;$\se {{\fh}^{\,g_1,g_2}} I: \se L I(F)\, \lra \, \se L I(F)$ as 
follows: for\, $h \in \se L I(F)$,
\vspace{-0.2cm}

\nhp \hfl $\se {{\fh}^{\,g_1,g_2}} I(h) = h\,\se g 1 \se g 2$. \hfl $\Box$
\edf
\vspace{-0.4cm}

\nhp \und{Note}. Since $Z(\se g i) \, \sub \; I = Z(h)\; (i = 1, 2)$, we have 
$Z(h\,\se g 1 \se g 2) = I$, whence $h\,\se g 1 \se g 2 \in \se L I$.
\vspace{-0.35cm}

\bfa \label{zeros-inv} With notation as in Definition $\ref{definvol}$, let 
$J \in$ {\em Spec}$(F)$ be such that $Z(\se g 1) \cup Z(\se g 2) \, \sub \, 
J \, \sub \; I$, and for \,$i = 1, 2,$ let \se {g'} i be the unique 
$\spez$-successor of\, \se g i  of level\, \y J. Then, $\se {{\fh}^{\,g_1,g_2}} I = 
\se {{\fh}^{\,g'_1,g'_2}} I$. Thus, in $\ref{definvol}$ we may assume 
$Z(\se g 1) = Z(\se g 2)$. 
\efa
\vspace{-0.4cm}

\nhp {\bf Proof.} Lemma \ref{prods} shows that $h\,\se g 1 \se g 2 = 
h\,\se {g'} 1 \se {g'} 2$, for $h \in \se L I$. \hfl $\Box$
\vspace{-0.35cm} 

\bth \label{propsinvol} With notation as in Definition $\ref{definvol}$, we have:
\vspace{-0.15cm}

\nhp $(a)$ $\se {{\fh}^{\,g_1,g_2}} I$ is an AOS-automorphism of\, \se L I.
\vspace{-0.15cm}

\nhp $(b)$ $\se {{\fh}^{\,g_1,g_2}} I$ is an involution: for $h \in \se L I$\,,\; 
$\se {{\fh}^{\,g_1,g_2}} I(\se {{\fh}^{\,g_1,g_2}} I(h)) = h$.
\vspace{-0.15cm}

\nhp $(c)$ For $i = 1, 2,$ let  \se h i be the unique $\spez$-successor of\, \se g i in  \se L I. 
Then,\; $\se {{\fh}^{\,g_1,g_2}} I(\se h 1) = \se h 2$.
\vspace{-0.15cm}

\nhp In particular,
\vspace{-0.15cm}

\nhp $(d)$ If \se g 1, \se g 2, have a common specialization \y h at some level 
$I \supseteq Z(\se g 1), Z(\se g 2)$, then\, \y h\, is a 
\vspace{-0.15cm}

\tab{\nhp $(d)$} fixed point of\, $\se {{\fh}^{\,g_1,g_2}} I$.
\vspace{-0.15cm}

\nhp $(e)$ Let $J \, \sub \, I$ be in $Spec(F)$. Assume $Z(\se g 1),Z(\se g 2) \, \sub \; J$, and let
$\se h 1 \in \se L {\!J},\; \se h 2 \in \se L I$. Then,
\vspace{-0.1cm}

\nhp \hfl $\se h 1 \spez \, \se h 2 \;\;\; \Ra \;\;\; \se {{\fh}^{\,g_1,g_2}} J(\se h 1) \spez \, \se {{\fh}^{\,g_1,g_2}} I(\se h 2)$. \hfl
\vspace{-0.15cm}

\eth
\vspace{-0.3cm}

For the proof of this Theorem we will need an improvement on \ref{levels}\,(a), 
valid for fans but not for arbitrary RSs; namely\,:
\vspace{-0.3cm}

\bfa \label{ImprovedInducedMap} Let \y F be a RS-fan, and \,\y I be an ideal of \y F. 
Any \,$g \in \se X {\!F}$ such that \und{$Z(g) \, \sub \, I$} induces a SG-character \,
$\h g: \se F I \, \lra \, \zt$, by setting \,$\h g \; \com \, \se {\pi} I = g$.  
\efa
\vspace{-0.45cm}

\nhp {\bf Proof.} The only delicate point is well-definedness: for\, 
$a \in F \setminus I$, \, $a \; \se {\sim} I \, 1 \;\, \Ra \;\, g(a) = 1.$ 
By \cite{DP5a}, Lemma 5.10,\: $a \; \se {\sim} I \, 1$ means $az = z$ for some 
$z \not\in \, I$ (see proof of \ref{RSG+levelMorph}); then $g(z) \neq 0$, and 
taking images under \y g in this equality yields $g(a) = 1$. \hfl $\Box$

\nhp {\bf Proof of Theorem \ref{propsinvol}.} We begin by proving:
\vspace{-0.15cm}

\nhp (b) For\, $h \in \se L I(F),\; \se {{\fh}^{\,g_1,g_2}} I(\se {{\fh}^{\,g_1,g_2}} I(h)) = 
h\,\se {g^2} 1 \, \se {g^2} 2$. But\; $h\,\se {g^2} 1 \, \se {g^2} 2 = h$; this is clear if 
$h(x) = 0$  $(x \in F)$; if $h(x) \neq 0$, then $\se g i(x) \neq 0$ (since $Z(\se g i) \, \sub \, Z(h)$), 
and hence $\se {g^2} i(x) = 1 \; (i = 1,2)$, proving the stated identity, and item (b).
\vspace{-0.15cm}

\noi (a)\: i) $\se {{\fh}^{\,g_1,g_2}} I$ is an AOS-morphism.
\vspace{-0.15cm}

\nhp Since \se F I is the RSG-fan dual to \se L I(\y F), we must show (see \ref{levels}\,(c)):
\vspace{-0.2cm}

\nhp (*) \; For every $\al \in \se F I$ there is $\bt \in \se F I$ such that 
$\h \al \; \com \; \se {{\fh}^{\,g_1,g_2}} I = \,\h \bt,$
\vspace{-0.15cm}

\nhp where $\h \al : \se X {F_I}\, \lra\, \zt$ denotes evaluation at \al. We claim that 
$\bt = \al\: \se {\h g} 1(\al) \, \se {\h g} 2(\al)$ does the job. By Fact
\ref{ImprovedInducedMap}, $\se {\h g} i(\al) \in \zt\; (i = 1,2)$, whence 
$\bt \in \se F I$. For $h \in \se L I(F)$ we have:
\vspace{-0.2cm}

\nhp \hfl $(\h \al \; \com \; \se {{\fh}^{\,g_1,g_2}} I)(h) = 
\h \al\,(h\,\se g 1 \se g 2) = h(\al)\, \se {\h g} 1(\al)\, \se {\h g} 2(\al) = 
h(\al\, \se {\h g} 1(\al)\, \se {\h g} 2(\al)) = h(\bt) = 
\h \bt\,(h)$, \hfl 
\vspace{-0.15cm}

\nhp as required. Note that (b) implies
\vspace{-0.15cm}

\nhp ii) \;$\se {{\fh}^{\,g_1,g_2}} I$ is bijective.
\vspace{-0.15cm}

\nhp iii) The dual map $(\se {{\fh}^{\,g_1,g_2}} I)^{*}: \se F I\, \lra \, \se F I$ 
is also bijective.
\vspace{-0.15cm}

\nhp Item (i) proves that, for $\al \in \se F I\,,\; (\se {{\fh}^{\,g_1,g_2}} I)^{*}(\al) = 
\al\: \se {\h g} 1(\al) \, \se {\h g} 2(\al)$. For injectivity, assume 
$\al\: \se {\h g} 1(\al) \, \se {\h g} 2(\al) = 1$; if $\se {\h g} 1(\al) \, 
\se {\h g} 2(\al) = -1$, then $\al = -1$, whence (as $\se {\h g} i$ is a 
SG-character), $\se {\h g} i(\al) = -1 \; (i = 1,2)$, and $\al = 1$, contradiction. 
Thus, $\se {\h g} 1(\al) \, \se {\h g} 2(\al) = 1$, which entails\,
$\al = 1$. For surjectivity, given $\bt \in \se F I$, set \,$\al = 
\bt\: \se {\h g} 1(\bt) \, \se {\h g} 2(\bt)$. Then, $\se {\h g} 1(\al) = 
\se {\h g} 2(\bt)$ and $\se {\h g} 2(\al) = \se {\h g} 1(\bt)$, whence 
$(\se {{\fh}^{\,g_1,g_2}} I)^{*}(\al) = \bt$.
\vspace{-0.15cm}

\nhp (c) We must prove $\se h 1\se g 1 \se g 2 = \se h 2$. This clearly\vspace{0.05cm} 
holds at any $x \in Z(\se h 1) = Z(\se h 2)$. If $x \not\in Z(\se h i)$\; 
$(i = 1,2)$, then $x \not\in Z(\se g i)$; since $\se g i \spez \, \se h i$, 
it follows $\se h i(x) = \se g i(x) \neq 0$ (Lemma \ref{char-specializ}\,(4)), 
and $\se h i(x)\, \se g i(x) = 1$; hence, $\se h 1\se g 1 \se g 2(x) = \se g 2(x) = 
\se h 2(x)$.
\vspace{-0.15cm}

\nhp (e) Lemma \ref{moreprods}\,(a) immediately implies\; $\se {h^2} 2 = 
\se h 2 \se h 1 \;\, \Ra \;\, (\se h 2 \se g 1 \se g 2)^2 = 
(\se h 2 \se g 1 \se g 2)(\se h 1 \se g 1 \se g 2)$. \hfl $\Box$

By use of these involutions we obtain a number of regularity results concerning 
the order structure of ARS-fans.
\vspace{-0.4cm}

\bpr \label{cardsetsS+C} Let \y F be a RS-fan. For $J \, \sub \, \se J 1 \, \sub\, 
\se J 2 \, \sub \, I$ in\, {\em Spec}$(F)$, and $h \in \se {S^I} {\!J}$ set:
\vspace{-0.2cm}

\nhp \hfl $B^{J_1,J_2}(h) = \{\, g \in \se {S^{J_2}} {J_1} \; \vert \; g \spez \, h\, \}$, \; 
and \; $A^{J_1,J_2}(h) = \{\, g \in \se {C^{J_2}} {J_1} \; \vert \; g \spez \, h\, \}$. \hfl
\vspace{-0.3cm}

\nhp Then,
\vspace{-0.2cm}
 
\nhp $(a)$ For \,$\se h 1, \se h 2 \in\, \se {S^I} {\!J}$\,, we have \, ${\mathrm {card}}\,(B^{J_1,J_2}(\se h 1)) 
= {\mathrm {card}}\,(B^{J_1,J_2}(\se h 2))$.
\vspace{-0.2cm}

\nhp $(b)$ For \,$\se h 1, \se h 2 \in\, \se {C^I} {\!J}$\,, we have \, ${\mathrm {card}}\,(A^{J_1,J_2}(\se h 1)) 
= {\mathrm {card}}\,(A^{J_1,J_2}(\se h 2))$.
\epr 
\vspace{-0.4cm}

\nhp {\bf Remark.} The assumptions of the Proposition guarantee that the sets $B^{J_1,J_2}(h)$ are non-empty. 
In fact, given $h \in \se {S^I} {\!J}$, there is $u\spez \, h$ so that $Z(u)\, \sub\, J$; set $J' = Z(u)$. 
Since  $J' \, \sub \, J \, \sub \, \se J 2$, \lbr
\y u has a unique $\spez$\,-\,successor \y g in $\se L {J_2}$. But $u \spez \, g,\, h$\, and\, 
$J_2 = Z(g) \, \sub \; I = Z(h)$ imply $g \spez \, h$ (Lemma \ref{char-zeroset}\,(3)). Since 
$J' \, \sub \, J \, \sub \, \se J 1$, we conclude that $g \in \se {S^{J_2}} {J_1}$, i.e., 
$g \in B^{J_1,J_2}(h)$. 
\vspace{-0.15cm}

The sets $A^{J_1,J_2}(h)$ may be empty for some choices of \,\y h\, and the {\se J i}'s. However, 
if $h \in\, \se {C^I} {\!J}$ and $\se J 1 = J$, we have $A^{J_1,J_2}(h) \neq \0$. Indeed, in this 
case the element $g \in \se {S^{J_2}} J$ constructed above is in $\se {C^{J_2}} J$, for if 
$g \in \se {S^{J_2}} {J'}$ for some $J' \subset J$, then $g \spez \, h$ would imply $h \in\, \se {S^I} {\!J'}$, 
contrary to the assumption $h \in\, \se {C^I} {\!J}$. \hfl $\Box$
\vspace{-0.1cm}

\nhp {\bf Proof of Proposition \ref{cardsetsS+C}.} (a) With \se J 1, \se J 2 as in the statement, 
write \se B i for $B^{J_1,J_2}(\se h i)\; (i = 1,2)$. The assumption $\se h i \in\, \se {S^I} {\!J}$ 
implies the existence of elements $\se u i \in \se X {\!F}$ so that $\se u i \spez \, \se h i$\, and \,
$Z(\se u i) \, \sub \: J$. Replacing \se u i by its unique successor of level \y J we may assume 
$Z(\se u i) = J$ (see \ref{zeros-inv}). We fix \se u i's with these properties throughout the proof, 
and for $J \, \sub \: J' \, \sub \: I$ we denote by \se {\fh} {\!\!J'} the involution 
$\se {{\fh}^{u_1,u_2}}{\!J'}$ of \se L {\!J'} defined in \ref{definvol}.
\vspace{-0.15cm}

Since the maps \se {\fh} {\!J'} are bijective, it is enough to prove $\se {\fh} {\!\!J_2}[\se B 1] = \se B 2$.
Further, since $\se {\fh}{\!\!J_2}$ is an involution it suffices just to prove the inclusion $\sub$, i.e.,
\vspace{-0.15cm}

\nhp (*) \hspace{0.5cm} $g \in \se {S^{J_2}} {\!J_1}$ \;and\; $g \spez \,\se h 1 \;\, \Ra \;\, 
\se {\fh}{\!\!J_2}(g) \spez \,\se h 2$ \;and\; $\se{\fh}{\!\!J_2}(g) \in \se {S^{J_2}} {\!J_1}$.
\vspace{-0.15cm}

\nhp (i)\; $\se {\fh}{\!\!J_2}(g) = g\,\se u 1 \se u 2 \spez \,\se h 2$.
\vspace{-0.15cm}

\nhp Immediate consequence of Lemma \ref{moreprods}\,(b), since $g, \se u 1 \spez \,\se h 1$ and 
$\se u 2 \spez \,\se h 2$.
\vspace{-0.2cm}

\nhp (ii) $\se{\fh}{\!\!J_2}(g) \in \, \se {S^{J_2}} {\!J_1}$.
\vspace{-0.2cm}

\nhp Since $g \in \, \se {S^{J_2}} {J_1}$, there is $v \spez \,g$ so that 
$Z(v) = \se J 1 \supseteq J = Z(\se u i)\; (i = 1,2)$; thus, \y v is in the domain 
of $\se{\fh}{\!Z(v)} = \, \se {{\fh}^{u_1,u_2}}{Z(v)}$, and Theorem \ref{propsinvol}\,(e)
gives $\se{\fh}{\!Z(v)}(v) \spez \, \se{\fh}{\!\!J_2}(g)$, proving (ii) and item (a). 
\vspace{-0.15cm}

\nhp (b) Write \se A i for $A^{J_1,J_2}(\se h i)\; (i = 1,2)$. As above, it suffices to prove 
the analogue of (*):
\vspace{-0.2cm}

\nhp (**) \hspace{0.5cm} $g \in \se {C^{J_2}} {J_1}$ \;and\; $g \spez \,\se h 1 \;\, \Ra \;\, 
\se {\fh}{\!\!J_2}(g) \spez \,\se h 2$ \;and\; $\se{\fh}{\!\!J_2}(g) \in \se {C^{J_2}} {\!J_1}$\,,
\vspace{-0.15cm}

\nhp where \se h 1, \se h 2 are now assumed to be in\, $\se {C^I} {\!J}$. In fact, by (*) it 
only remains to show:
\vspace{-0.15cm}

\nhp (iii) There is no \,$w \in\, \se X {\!F}$ such that \,$Z(w) \subset \se J 1$\, and\, $w \spez \, 
\se {\fh}{\!\!J_2}(g)$.
\vspace{-0.15cm}

\nhp Otherwise, we would have $w \spez \, \se {\fh}{\!\!J_2}(g) \spez \, \se h 2$ (the last relation holding 
by (*)). Since $\se h 2 \in \, \se {C^I} {\!J}$, we get $J \, \sub \, Z(w)$, and since $Z(\se u i) = J$, 
$\se{\fh}{\!Z(w)}(w)$ is defined. Theorem \ref{propsinvol}\,(e) applied to the first of the preceding 
inequalities yields: $\se{\fh}{\!Z(w)}(w) \spez \se {\fh}{\!\!J_2}(\se {\fh}{\!\!J_2}(g)) = g$. This, 
together with $\se{\fh}{\!Z(w)}(w) \in \se L {Z(w)}$\, and $Z(w) \subset \se J 1$, contradicts the 
assumption $g \in \se {C^{J_2}} {\!J_1}$, proving (iii), and item (b). \hfl $\Box$

A slight variant of the argument proving Proposition \ref{cardsetsS+C} yields:
\vspace{-0.4cm}

\bpr \label{permsetsS+C} Let \y F be a RS-fan and let $J \, \sub \, I$ be in 
$\Spec(F)$. For \,$\se g 1, \se g 2 \in \se X {\!F}$ such that \,$Z(\se g i) \, 
\sub \, J$ \, $(i = 1,2)$, the map \,$\se {{\fh}^{\,g_1,g_2}}I$ is a permutation 
of \,$\se {S^I} {\!J}$\, and of \,$\se {C^I} {\!J}$. \hfl $\Box$
\epr
\vspace{-0.3cm}

For a RS-fan, \y F, and $h \in \se X {\!F}$, we denote by \se P h = 
$\{\,g \in \se X {\!F} \, \vert \; g \spez\, h \,\}$ the root-system of 
predecessors of \y h under specialization. We begin by proving:
\vspace{-0.4cm}

\bpr \label{fan-predec} $(1)$ \se P h is an ARS-fan. In particular,
\vspace{-0.15cm}

\nhp $(2)$ Any connected component of an ARS-fan is an ARS-fan.
\epr

\vspace{-0.35cm}

\nhp {\bf Proof.} (1) Lemma \ref{char-specializ}\,(2) shows that\, $g \spez\, h$\, iff\, 
$T = h^{-1}[1] \, \sub \; g^{-1}[1]$. With notation as in \cite{M}, \S\,6.3, p. 110, and   
\S\,6.6, p. 126, the latter condition just means $g \in U(T)$, i.e., \se P h is the saturated 
set $U(T) (= W(T) \, \cap \, U(T^2))$. \cite{M}, Cor. 6.6.8, p. 126, proves that sets of this 
form are ARSs. Lemma \ref{moreprods}\,(a) shows that it is closed under products of three 
elements, hence a fan by the results of \cite{DP5a}, \S\,3.
\vspace{-0.1cm}

\noi (2) Follows from (1) by taking $h$ to be the (unique) $\spez$-top element of the given 
connected component. \hfl $\Box$

Continuing the analysis of (ARS-)fans of the form \se P h, we show:

\vspace{-0.4cm}

\bth \label{embed-predec} Let \y F be a RS-fan and let $J \, \sub \, I$ be in\, {\em Spec}$(F)$. 
Let $\se h 1 \in \se {C^I} {\!J},\; \se h 2 \in \se {S^I} {\!J}$. For $i = 1,2,$ we write 
\se P i for \se P {h_i}. Then,
\vspace{-0.15cm}

\nhp $(1)$ There is an ARS-embedding \,\fh\: of\, \se P 1 into \se P 2. Further, $\fh[\se P 1] = 
\{\,u \in \se P 2 \, \vert \; J \, \sub \; Z(u)\:\}$. In\vspace{0.05cm} particular, \fh\ is an 
order-embedding of\, $(\se P 1, \spez)$ into $(\se P 2, \spez)$.
\vspace{-0.15cm}

\nhp $(2)$ If, in addition, $\se h 2 \in \se {C^I} {\!J}$, then\, \fh\ is an isomorphism of ARSs.
\eth

\vspace{-0.3cm}

\nhp {\bf Proof.} Our assumption on the \se h i's guarantees the existence of $\se u 1, \se u 2 \in \se L {\!J}$ 
so that $\se u i \spez\, \se h i\; (i = 1,2)$. For $J \, \sub \,J'\, \sub\, I$ in $\textrm{Spec}(F)$ 
let \se {\fh} {\!\!J'} denote the involution \se {\fh^{u_1,u_2}} {\!\!J'} of \se L {J'} (Definition \ref{definvol}).
\vspace{-0.15cm}

\nhp (1) We construct $\fh : \se P 1 \, \lra \, \se P 2$ by ``collecting together" all the relevant maps 
\se {\fh} {\!\!J'}\, $(J \, \sub \,J'\, \sub\, I)$: given $g \in \se L {J'}\,,\, g \spez\, \se h 1$, we set

\vspace{-0.25cm}

\nhp \hfl $\fh(g) = \se {\fh} {\!\!J'}(g)$. \hfl

\vspace{-0.15cm}

\nhp Since the levels \se L {\!J'} are pairwise disjoint,\, \fh\ is well-defined.
\vspace{-0.15cm}

\nhp i) \; $\fh[\se P 1]\, \sub \, \se P 2$.
\vspace{-0.15cm}

\nhp By Theorem \ref{propsinvol}\,(e),\, $g \spez\, \se h 1$ implies $\se {\fh} {\!\!J'}(g) \spez\, 
\se {\fh} {\!I}(\se h 1)$. Since \se h i is the unique successor of \se u i at level \y I, 
\ref{propsinvol}\,(c) yields $\se {\fh} {\!\!I}(\se h 1) = \se h 2$, whence 
$\se {\fh} {\!\!J'}(g) \spez\, \se h 2$, as required. Note this also gives $J \, \sub \; J' = Z(\fh(g))$.
\vspace{-0.15cm}

\nhp ii)\, $\{\,u \in \se P 2 \, \vert \; J \, \sub \; Z(u)\:\} \, \sub \; \fh[\se P 1]$.
\vspace{-0.15cm}

\nhp Let \y u be in the left-hand side, with $J' = Z(u)$, say. Set $v = \se {\fh} {\!\!J'}(u)$; then,
$\se {\fh} {\!\!J'}(v) = u$ (\ref{propsinvol}\,(b)). By \ref{moreprods}\,(b), $\se u 1 \spez\; \se h 1$ 
and\vspace{0.05cm} $u, \se u 2 \spez\, \se h 2$ imply $u\,\se u 1\,\se u 2 = \se {\fh} {\!\!J'}(u) = v \spez\, \se h 1$,
i.e., $v \in \se P 1$. Hence $\fh(v) = u \in\, \fh[\se P 1]$.
\vspace{-0.15cm}

\nhp iii) \fh\, is injective.
\vspace{-0.15cm}

\nhp This is clear using \ref{propsinvol}\,(b), since $Z(\fh(g)) = Z(g)$ for $g \in \se P 1$.
\vspace{-0.15cm}

\nhp iv)\, \fh\, is an ARS-morphism.
\vspace{-0.15cm}

\nhp The proof is similar to that of item (a) in Theorem \ref{propsinvol}. The statement to be proved is:
\vspace{-0.2cm}

\nhp (\dag) For every $a \in F$ there is $b \in F$ such that $(\h{a/\se T 2}) \,\com\, \fh = 
\h{b/\se T 1}$,
\vspace{-0.2cm}

\nhp where, for $i = 1,2$, $\se T i = \se h i^{-1}[1],\: \se P i = U(\se T i),\:  
\h{a/\se T 2}: \se P 2 \, \lra \;{\bf 3}$\, is the evaluation map: 
$\h{a/\se T 2}\,(g) = \h g\,(a/\se T 2) = g(a)$, for $g \in \se P 2$, and similarly 
for $\h{b/\se T 1}: \se P 1 \, \lra \;{\bf 3}$. (Note that $g \in \se P 2 = U(\se T 2)$ 
ensures that $\h{a/\se T 2}$ depends only on the congruence class of \y a modulo \se T 2.)
\vspace{-0.15cm}

The conclusion of (\dag) can equivalently be written as\, $\h {\fh(g)}(a/\se T 2) = \h g\,(b/\se T 1)$, 
i.e.,\, $(\se u 1\se u 2\,g)(a)$ $= g(b)$. Since $\se u i(a) \in \{0,1,-1\}\; (i=1,2)$, it is clear 
that the element $b = a\,\se u 1(a)\se u 2(a) \in F$ verifies (\dag); see \ref{propsinvol}\,(a).
\vspace{-0.15cm}

\nhp (2) Since $\se h 2 \in \se {C^I} {\!J}$, the preceding construction can be performed with the roles of 
\se h 1 and \se h 2 reversed. Routine verification using \ref{propsinvol}(b) shows that the map obtained 
is\, $\fh^{-1}$, which then is an ARS-morphism, proving that \,\fh\ is an ARS-isomorphism. \hfl $\Box$

Proposition \ref{cardsetsS+C} and Theorem \ref{embed-predec} provide 
significant information on the structure of the connected components of ARS-fans.

\vspace{-0.45cm}

\bdfr \label{components} {\em (a) Let $(X,\preceq\,)$ be a root-system and let 
$\se g 1, \se g 2 \in X$. Define:
\vspace{-0.15cm}

\nhp \hfl $\se g 1\, \se {\equiv} C\, \se g 2$ \; iff \; \se g 1, \se g 2 have
a common $\preceq$\,-\,upper bound. \hfl
\vspace{-0.15cm}

\nhp $\se {\equiv} C\,$ is an equivalence relation; its classes are called 
{\bf connected components} of\, $(X,\preceq\,)$.
\vspace{-0.15cm}

\noi (b) The $\spez$\,-\,top elements of the connected components of an 
ARS-fan $(X,F)$ have all the same level, namely the level determined by the 
 maximal ideal \y M of \y F; cf. Proposition \ref{charfan:(1)implies(2.ii)}\,(3). 

\vspace{-0.15cm}

\nhp (c) Since every connected component of an ARS-fan is itself an ARS-fan,
\ref{fan-predec}\,(2), the zero-sets of its elements attain a lowest level, which 
can be explicitly determined, cf. Proposition \ref{componlowestlevel} below. However, 
different components may have different lowest levels, see Corollary \ref{embedcompon}. }
\hfl $\Box$
\edfr
\vspace{-0.3cm}

\nhp {\bf Notation.} The sets \se L I, $\se {S^I} {\!J}$ and $\se {C^I} {\!J}$ defined in 
\ref{depthlev} and \ref{setsS+C} relativize in an obvious way to the connected 
components of a fan $(X,F)$; if \y K is such a component and $J \, \sub \, I$ 
are in Spec(\y F) we set:
\vspace{-0.4cm}

\nhp \hfl $\se L I(K) = \se L I \, \cap \, K$, \;\;\; $\se {S^I} {\!J}(K) = \se {S^I} {\!J} \, 
\cap \, K$,\;\;\; and\;\;\; $\se {C^I} {\!J}(K) = \se {C^I} {\!J} \, \cap \, K$. \hfl
\vspace{-0.2cm}

\nhp Note that some (or all) of these sets may be empty, depending on 
\y I, \y J and the component \y K. Clearly, if \se h 0 is the 
$\spez$\,-\,top element of \y K, we have $\se L I(K) = 
\{\, g \in \se L I \, \vert \, g \spez \,\se h 0 \,\}$, and similarly 
for $\se {S^I} {\!J}(K)$\, and\, $\se {C^I} {\!J}(K)$. $\se L I(K) \neq \0$ 
just means that \y K ``reaches at least"  the \y I-th level of \y X 
(possibly lower). \hfl $\Box$

\vspace{-0.35cm}

\bpr \label{componlowestlevel} Let \y K be a connected component of an 
ARS-fan $(X,F)$. Let \se h 0 be the $\spez$\,-\,top element of \y K, and 
let $T = \se {h^{-1}} 0[1]$. Then, the lowest\vspace*{0.05cm} level of \y K 
$($i.e., the smallest ideal \y I of \y F such that $\se L I(K) \neq \0)$ is\, 
$I = \Ga \, \cap \, -\Ga$, where \Ga\ is the saturated subsemigroup of\, 
\y F generated by\, ${\mathrm {Id}}(F) \cdot T$.
\epr
\vspace{-0.4cm}

\nhp \und{Note}. The subsemigroup Id$(F) \cdot T$ may not be saturated, 
since ${\mathrm {Id}}(F) \cdot T \, \cap \, -({\mathrm {Id}}(F) \cdot T)$ 
is not, in general, an ideal; see \cite{DP5a}, Cor. 3.10\,(2).
\vspace{-0.15cm}

\nhp {\bf Proof.} With notation as in \ref{fan-predec}, we have $K = \se P {h_0} = U(T) = 
\{ g \in X \, \vert \, g\, \lceil\, T = 1 \} =$ the ARS $\se X {F/T}$ (where 
$F/T = F/{\se {\sim} K}$, with $\se {\sim} K$ denoting the congruence on \y F induced by \y K). 
Let $\se {\pi} T : F \, \lra \ F/T$ be the quotient map. The lowest level of $\se X {F/T}$ is $\{0\}$; 
with \y K identified to a subset of \y X via the map $g \, \mapsto \h g \;\; (\h g \; \com \, \se {\pi} T = g)$, 
the corresponding ideal of \y F is  $\se {\pi^{-1}} T[0] = 
\{ a \in F \, \vert \, a\, \se {\sim} K \, 0 \}$. Then, with the ideal $I$ defined in the statement, 
we must\vspace{0.08cm} prove, for $a \in F$\,:
\vspace{-0.5cm}

\nhp \hfl $a\; \se {\sim} K \, 0 \;\; \Longleftrightarrow \;\; a \in I$. \hfl
\vspace{-0.15cm}

\nhp ($\La$) This follows from $I \, \sub \, Z(g)$ for all $g \in K$. Since $g\, \lceil\, T = 1$,
we get ${\mathrm {Id}} \cdot T \; \sub \; P(g) = g^{-1}[0,1]$; since $P(g)$ is a saturated 
subsemigroup, it comes $\Ga \, \sub \, P(g)$. Hence, $x \in I = \Ga \, \cap \, -\Ga$ implies $g(x) = 0$.
\vspace{-0.15cm}

\nhp ($\Ra$) Assume $a \not\in I$. In order to prove $a\; \se {\not\sim} K \, 0$ we construct 
a character $g \in X$ such that $g\, \lceil\, T = 1$ and $g(a) \neq 0$ (i.e., $g(a^2) = 1$). 
The ideal \y I is prime and saturated (\cite{DP5a}, 3.10\,(1)). Since $I = \Ga \, \cap \, -\Ga$, 
there is a saturated subsemigroup $S$ of $F$ such that $\Ga\, \sub\, S$ and $S$ maximal with 
$S\, \cap\, -S = I$. By \cite{DP1}, Lemma 3.5, p. 114, $S\, \cup\, -S = F$, and $S$ defines a 
character $g \in X$ with $Z(g) = I$, by setting 
$g\, \lceil\, (S \setminus -S) = 1,\; g\, \lceil\, (-S \setminus S) = -1$ and $g\, \lceil\, I = 0$.
Note that we have,

\vspace{-0.15cm}

\nhp ($\dag$) \; $I\, \cap\, a^2T = \0$.
\vspace{-0.2cm}

\noi Otherwise, there is $t \in T$ such that $a^2t \in I$; since $I$ is prime and $a \not\in I$, we get
$t \in I$, contradicting $T\, \cap\, I = \se {h^{-1}} 0[1]\, \cap\, Z(\se h 0) = \0$.
\vspace{-0.2cm}

Since $a^2T\, \sub\, S$, (\dag) implies $-S\, \cap\, a^2T = \0$, whence $g\, \lceil\, a^2T = 1$ by the
definition of $g$. \hfl $\Box$

Proposition \ref{cardsetsS+C} implies:

\vspace{-0.4cm}

\bco \label{cardcomponents} Let \,$(X,F)$ be an ARS-fan and let \se K 1, \se K 2 be connected 
components of\, $(X,F)$. Then,
\vspace{-0.15cm}

\nhp $(1)$ Let $I \in$ {\em Spec}$(F)$; if\, $\se L I(\se K i) \neq \0$ for\, $i = 1, 2$, then\; 
${\mathrm {card}}\, (\se L I(\se K 1)) = {\mathrm {card}}\,(\se L I(\se K 2))$. 
\vspace{-0.15cm}

\nhp $(2)$ Let $J \, \sub \; J'$ be in\, {\em Spec}$(F)$, and assume $\se L {\!J}(\se K i) \neq \0\;\; (i = 1,2)$. 
Then,\, ${\mathrm {card}}\,(\se {S^{J'}} J(\se K 1)) = {\mathrm {card}}\,(\se {S^{J'}} J(\se K 2))$. 
\eco

\vspace{-0.4cm}

\nhp {\bf Proof.} (1) follows from (2), as $\se L I = \se {S^I} I$.
\vspace{-0.15cm}

\nhp (2) Fix $i \in \{1,2\}$. Let \se h i be the $\spez$\,-\,top element of \se K i. The assumption 
$\se L {\!J}(\se K i) \neq \0$ implies that the sets $\se {S^{J'}} J(\se K i) = 
\{\, g \in \se {S^{J'}} J \, \vert \, g \spez \,\se h i \,\}$ are non-empty. Now, applying Proposition 
\ref{cardsetsS+C}(a) with \y I = \y M (= the maximal ideal of \y F), $\se J 1 = J$, $\se J 2 = J'$ 
we have $B^{J,J'}(\se h i) = \{\, g \in \se {S^{J'}} J \; \vert \; g \spez \,\se h i\, \} = 
\se {S^{J'}} J(\se K i)$, and the result follows. \hfl $\Box$
\vspace{-0.1cm}

\nhp {\bf Remark.} Assertion (2) of Corollary \ref{cardcomponents} fails, in general, if the sets 
$\se {S^{J'}} J(\se K i)$ are replaced by $\se {C^{J'}} J(\se K i)$, even if both sets 
$\se {C^{J'}} J(\se K i),\; i = 1, 2,$\, are assumed non-empty. The snag is that 
$\se {C^{J'}} J(\se K i) \neq \0$\, \und{does not} imply that the $\spez$\,-\,top element 
\se h i of \se K i belongs to $\se {C^M} J(\se K i)$, a condition required for Proposition 
\ref{cardsetsS+C}(b) to apply. It is easy to construct counterexamples. \hfl $\Box$


Theorem \ref{embed-predec} gives:

\vspace{-0.35cm}

\bco \label{embedcompon} Let \se K 1, \se K 2 be connected components of the ARS-fan\, $(X,F)$. Let 
$\se I 1, \se I 2 \in$ {\em Spec}$(F)$ be the lowest levels of \se K 1, \se K 2, resp. $($cf.  
$\ref{componlowestlevel}\,)$. Then,
\vspace{-0.15cm}

\nhp $(1)$ If \,$\se I 2 \,\sub\; \se I 1$, then\, \se K 1 endowed with the specialization order is $($order-$)$isomorphic to the root-system obtained by deleting all levels $I \subset \se I 1$ in \se K 2.
\vspace{-0.15cm}

\nhp $(2)$ If\, $\se I 1 = \se I 2$, then\, \se K 1, \se K 2 are order-isomorphic.
\eco

\vspace{-0.35cm}

\nhp {\bf Proof.} (1) Use Theorem \ref{embed-predec}\,(1) with \y I = \y M = the maximal ideal of 
\y F,\, \y J = \se I 1,\, and \se h 1, \se h 2 the $\spez$\,-\,top elements of \se K 1, \se K 2, 
resp. The ARS-embedding $\fh : \se K 1\, \lra \, \se K 2$ constructed therein verifies 
$\fh[\se K 1] = \{\,u \in \se K 2 \, \vert \; \se I 1 \, \sub \; Z(u)\:\}$, which is exactly 
statement (1).
\vspace{-0.15cm}

\nhp (2) follows from Theorem \ref{embed-predec}\,(2). \hfl $\Box$
\vspace{-0.5cm}

\bct \label{impossibles} {\bf Some impossible configurations.}

The preceding results show that there are strong constraints on the order structure of ARS-fans, 
especially when there is more than one connected component. We include a few examples to help 
the reader visualize the extent of those restrictions.
\vspace{-0.15cm}

\nhp (1) A configuration like

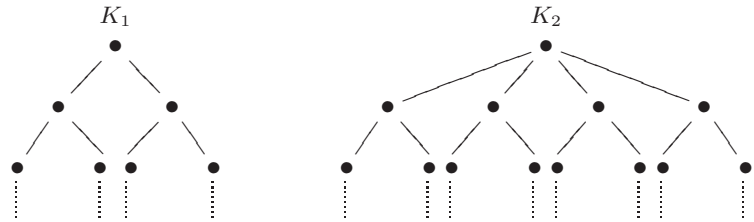
\vspace{-1cm}

$$\xymatrix@R=12pt@C=4pt{
  & & & & \\
  & & {\bullet} \ulab{$K_1$} \ar@{-}[dl] \ar@{-}[dr] \\
  & {\bullet} \ar@{-}[dl] \ar@{-}[dr] & & {\bullet} \ar@{-}[dl] \ar@{-}[dr] \\
  {\bullet} \ar@{.}@<-0.5pt>[d] & & {\bullet} \ar@{.}@<-6pt>[d] \;\; {\bullet} \ar@{.}@<+4pt>[d] & & {\bullet} \ar@{.}@<-0.5pt>[d]\\ 
  & & & & & \\ } \;\;\;\;\;\;
\xymatrix@R=12pt@C=4pt{
& & & & & & & \\
& & & & & {\bullet} \ulab{$K_2$} \ar@{-}[dlll] \ar@{-}[dl] \ar@{-}[dr] \ar@{-}[drrr] \\
& & {\bullet} \ar@{-}[dl] \ar@{-}[dr] & & {\bullet} \ar@{-}[dl] \ar@{-}[dr] & & {\bullet} \ar@{-}[dl] \ar@{-}[dr] & & {\bullet} \ar@{-}[dl] \ar@{-}[dr] \\ 
& {\bullet} \ar@{.}@<-0.7pt>[d] & & {\bullet} \; {\bullet} \ar@{.}@<-5pt>[d] \ar@{.}@<+3.5pt>[d] & & {\bullet} \; {\bullet} \ar@{.}@<-5pt>[d] \ar@{.}@<+3.5pt>[d] & & {\bullet} \; {\bullet} \ar@{.}@<-5pt>[d] \ar@{.}@<+3.5pt>[d] & & {\bullet} \ar@{.}@<-1pt>[d]\\ 
& & & & & & & & & \\ }
$$

\vspace{-0.5cm}

\nhp contradicts Corollary \ref{cardcomponents}\,(1).
\vspace{-0.2cm}

\nhp (2) The four-component configuration

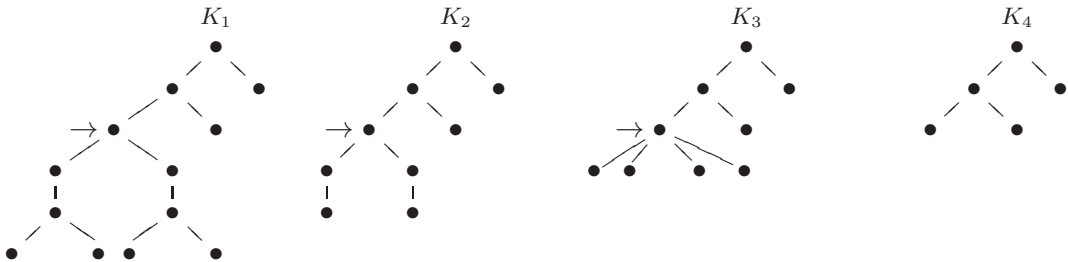
\vspace{-0.95cm}

$$\xymatrix@=1ex{
& & & & & & \\
& & & & {\bullet} \ulab{$K_1$} \ar@{-}[dl] \ar@{-}[dr] \\
& & & {\bullet} \ar@{-}[dl] \ar@{-}[dr] & & {\bullet} \\
& & {\hspace*{-0.5cm}\rightarrow \bullet} \ar@{-}[dl] \ar@{-}[dr] & & {\bullet} \\
& {\bullet} \ar@{-}[d] & & {\bullet}\ar@{-}[d] \\
& {\bullet} \ar@{-}[dl] \ar@{-}[dr] & & {\bullet} \ar@{-}[dl] \ar@{-}[dr] \\
{\bullet} & & {\bullet} \;\; {\bullet} & & {\bullet} \\ 
& & & & & \\ } \;
\xymatrix@=1ex{
& & & & & \\
& & & {\bullet} \ulab{$K_2$} \ar@{-}[dl] \ar@{-}[dr] \\
& & {\bullet} \ar@{-}[dl] \ar@{-}[dr] & & {\bullet} \\
& {\hspace*{-0.5cm}\rightarrow \bullet} \ar@{-}[dl] \ar@{-}[dr] & & {\bullet} \\
{\bullet} \ar@{-}[d] & & {\bullet}\ar@{-}[d] \\
{\bullet}  & & {\bullet} \\
& & & & \\ }\;\;
\xymatrix@=1ex{
& & & & & & \\
& & & & & {\bullet} \ulab{$K_3$} \ar@{-}[dl] \ar@{-}[dr] \\
& & &  & {\bullet} \ar@{-}[dl] \ar@{-}[dr] & & {\bullet} & & & \\
& & & {\hspace*{-0.5cm}\rightarrow \bullet} \ar@{-}[dll]  \ar@{-}[dl] \ar@{-}[dr] \ar@{-}[drr] & & {\bullet} \\
& & {\hspace*{-0.8cm} \bullet} & {\hspace*{-0.8cm} \bullet} & {\hspace*{-0.1cm} \bullet} & {\hspace*{-0.05cm} \bullet} \\
& & & &  \\ } \;\;\; 
\xymatrix@=1ex{
& & & & \\
& & {\bullet} \ulab{$K_4$} \ar@{-}[dl] \ar@{-}[dr] \\
& {\bullet} \ar@{-}[dl] \ar@{-}[dr] & & {\bullet} \\
{\bullet} & & {\bullet} \\
& & & \\ } 
$$

\vspace{-0.85cm}

\nhp (where the components pairwise verify the conclusion of \ref{cardcomponents}\,(2)) is also 
impossible:\, card\,$(\se {S^{\,3}} 4)$ $= 3$ is not a power of 2, and hence $\se {S^{\,3}} 4$ 
(shown with arrows) cannot be an AOS-fan (see Corollary \ref{fansS}). However, the same 
configuration with \se K 3 replaced by another copy of \se K 4 does not clash with either 
\ref{cardcomponents} or \ref{embedcompon}.
\vspace{-0.15cm}

\nhp \und{Note}. Our notation here (and below) follows the convention introduced in \ref{finfannotat} 
for finite fans. Thus, $\se {S^{\,3}} 4$ stands for the set $\se {S^{\,I_3}} {\,I_4}$, see 
\ref{setsS+C} and \ref{notatff}.
\vspace{-0.1cm}

\nhp (3) The two-component root-system

\vspace{-1.3cm} 

\nhp \hspace*{1cm} \parbox[t]{230pt} {
$$\xymatrix@=1ex{
& & & & & & & & & & & & & & & & & & & & & & & & & & & & \\
& & & & & & & & & & & & & & & & {\bullet} \ulab{$K_1$} \ar@{-}[ddlllllllllllll] \ar@{-}[ddllllll] \ar@{-}[ddr] \ar@{-}[ddrrrr] & & & & & & & & & \\
& & & & & & & & & & & & & & & & & & & & & & & & & & \\
& & & {\bullet} \ar@{-}[ddll] \ar@{-}[ddrr] & & & & & & & {\bullet} \ar@{-}[ddll] \ar@{-}[ddrr] & & & & & & & {\bullet} & & & {\bullet} & & & & &  \\
& & & & & & & & & & & & & & & & & & & & & & & & & \\
& {\bullet} \ar@{-}[d] & & & & {\bullet} \ar@{-}[d] & & & {\bullet} \ar@{-}[d] & & & & {\bullet} \ar@{-}[d] & & & & & & & & & & & & &  \\
& {\bullet} \ar@{-}[ddl] \ar@{-}[ddr] & & & & {\bullet} \ar@{-}[ddl] \ar@{-}[ddr] & & & {\bullet} \ar@{-}[ddl] \ar@{-}[ddr] & & & & {\bullet} \ar@{-}[ddl] \ar@{-}[ddr] & & & & & & & & & & & & & \\
& & & & & & & & & & & & & & & & & & & & & & & & & & & & & & \\
{\bullet} & & {\bullet} & & {\bullet} & & {\bullet} & {\bullet} & & {\bullet} & & {\bullet} & & {\bullet} & & & & & & & & & & & &  \\ }
$$
}

\hspace{4.5cm} \parbox[t]{20pt}{\vspace*{1.1cm} 

1 \\ \vspace*{0.12cm} 

2 \\ \vspace*{-0.17cm}

3 \\ \vspace*{-0.3cm}

4 \\ \vspace*{-0.15cm}

5}

\vspace{-1cm} 

\nhp \hspace*{1cm} \parbox[t]{230pt} {
$$\xymatrix@=1ex{
& & & & & & & & & & & & & & & & & & \\
& & & & & & & & & & & & & & & & {\bullet} \ulab{$K_2$} \ar@{-}[ddlllllllllllll] \ar@{-}[ddllllll] \ar@{-}[ddr] \ar@{-}[ddrrrr] & & & & & & & \\
& & & & & & & & & & & & & & & & & & & & & & & & \\
& & & {\bullet} \ar@{-}[ddll] \ar@{-}[ddrr] & & & & & & & {\bullet} \ar@{-}[ddll] \ar@{-}[ddrr] & & & & & & & {\bullet} & & & {\bullet} \\
& & & & & & & & & & & & & & & & & & & & & & & & \\
& {\bullet} \ar@{-}[d] & & & & {\bullet} \ar@{-}[d] & & & {\bullet} & & & & {\bullet} & & & & & & & & \\
& {\bullet} \ar@{-}[ddl] \ar@{-}[ddr] & & & & {\bullet} \ar@{-}[ddl] \ar@{-}[ddr] & & & & & & & & & & & & & & & & & \\
& & & & & & & & & & & & & & & & & & & & & & & & \\
{\bullet} & & {\bullet} & & {\bullet} & & {\bullet} & & & & & & & & & & & & & & \\ } 
$$
}
\hspace*{5.2cm} \parbox[t]{20pt} {
\vspace{-0.4cm}
$$\hspace{-2.3cm}\xymatrix@=1ex{
& & \\
& & \\
& & 1 \\
& & & \\
& & 2 \\
& & & \\
& & 3 \\
& & 4 \\
& & 5 \\ }
$$
}

\vspace{-0.4cm}
\nhp contradicts both Corollary \ref{cardcomponents}\,(2) (card\,$(\se {S^{\,3}} {\,4}(\se K 1)) = 4$, but\: 
card\,$(\se {S^{\,3}} {\,4} (\se K 2)) = 2$) and Corollary \ref{embedcompon} (\se K 1 and \se K 2 have the 
same ``length" but are not order-isomorphic).  \hfl $\Box$
\ect

\vspace{-0.7cm}

\section{The specialization root-system of finite ARS-fans}\label{specializffans}

\vspace{-0.3cm}

In this section we shall mostly deal with finite fans in the categories {\bf ARS} and
{\bf RS}. Our main result is Theorem \ref{isothffans} ---the isomorphism theorem
for finite ARS-fans--- which proves that, in this case, the order of
specialization alone determines the isomorphism type. The proof depends on the
notion of a ``standard generating system" which we introduce in \ref{stgen}. \hfl $\Box$

\vspace{-0.5cm}

\bct {\bf Notation and Reminder} \label{notatff} (a) Notation \ref{finfannotat} for finite (ARS- 
and RS-)fans is used systematically in this section, adapted in a self-explanatory way; e.g.,
for $1 \leq k \leq j \leq n = \ell(\se X F)$,
\se L k (or \se L k(\se X F), if necessary), will stand for the level $\se L
{I_k},\; \se {S^k} j$ for the set $\se {S^{I_k}} {I_j}$, etc.
\vspace{-0.15cm}

\nhp (b) Recall that the AOSs have a combinatorial geometric (matroid) structure; it was
introduced in \cite{D1} and \cite{D2} for spaces of orders of fields, and later
generalized to abstract order spaces in \cite{Li}. In general, ARSs \und{do not
possess such a structure}. Thus, combinatorial geometric notions such as {\it
dependent set, independent set, basis, closed set, closure, dimension}, etc.,
will always refer to the above-mentioned combinatorial geometric
structure, and apply \und{only} to AOSs. For the definition and the mutual
relationships, in the general context of matroid theory, of combinatorial
notions such as those just mentioned, the reader is referred to \cite{Wh}.
\vspace{-0.15cm}

Since the combinatorial geometric structure of any AOS is isomorphic to that of
a set of vectors in a (possibly infinite-dimensional) vector space over the
two-element field\, \Ft\, with the structure induced by linear dependence (cf.
\cite{D1}, Thm. 3.1, p. 618), the notions above coincide with the corresponding
notions over vector spaces. For example, a subset $A\, \sub\, X$ of an AOS
$(X,G,-1)$ (\y G a group of exponent 2) is \und{dependent} iff there are
pairwise distinct elements $g, \se g 1, \dots, \se g r \in A\, (r \geq 2)$, such
that\, $g = \se g 1 \cdot\ \dots\ \cdot \se g r$ (as characters of \y G).
Since functions in \y X send $-1$ to $-1$, this functional identity can only hold 
if \y r is odd. Likewise, \y A is \und{closed} iff the product of any odd number 
of members of \y A belongs to \y A. \hfl $\Box$
\vspace{-0.15cm}

\noi {\bf Warning.} In this section the words {\it closed set} and {\it closure} are
used \und{only} in the combinatorial geometric sense just defined. \hfl $\Box$

\ect

\vspace{-0.8cm}

\ble \label{dep-fans} Let $(X,F)$ be an ARS-fan $($not necessarily finite$)$.
Let $J\, \sub\, I$ be in\, {\em Spec(}\y F$)$, and let $A\, \sub\, \se L J$,
$B\, \sub\, \se L I$, be sets such that:
\vspace{-0.15cm}

\nhp $($i$)$\, The unique $\spez$-\,successor in $\se L I$ of each $g \in A$
belongs to \y B.
\vspace{-0.15cm}

\nhp $($ii$)$ Every\, $h \in B$\, has a unique $\spez$-\,predecessor in \y A.
\vspace{-0.15cm}

\nhp Then, \y A dependent\; $\Ra$\; \y B dependent.
\ele

\vspace{-0.4cm}
\nhp {\bf Proof.} By assumption there are pairwise distinct elements $g,\, \se g
1, \dots,\, \se g r \in A$ such that\, $g = \se g 1 \cdot\ \dots\ \cdot \se g
r$; as observed above, \y r is odd $\geq 3$. Let $h,\, \se h 1, \dots,\, \se h
r$ be the unique successors of $g,\, \se g 1, \dots,\, \se g r$, resp., in \y B
coming from (i); thus,\, $g \spez\, h$ and $\se g i \spez\, \se h i$, for $i
= 1, \dots, r$. By \ref{moreprods}(a) we have $g = \se g 1 \cdot\ \dots\ \cdot\,
\se g r \spez\, \se h 1 \cdot\ \dots\ \cdot\, \se h r$. Since $\se h 1 \cdot\
\dots\ \cdot \se h r \in \se L I$ (\y r is odd) and \y g has a unique
$\spez$\,-\,successor in \se L I, we get $h = \se h 1 \cdot\ \dots\ \cdot \se
h r$.

By assumption (ii), every element in \y A is the {\it unique} predecessor of an
element in \y B. Since $\se g i \neq \se g j$, we get $\se h i \neq \se h j$ for
$1 \leq i \neq j \leq r$; likewise, $h \neq \se h i$ for $i = 1, \dots, r$. This
proves that \y h is the product of \y r {\it distinct} elements in \y B, and hence
that \y B is dependent. \hfl $\Box$
\vspace{-0.4cm}

\bpr \label{choicebas} {\em (Choice of basis).} Let $(X,F)$ be a finite ARS-fan;
let $1 \leq k < n = \ell(X)$. Let\, \cG\, be an arbitrary AOS-subfan of 
$\se L {k + 1} = \se L {k + 1}(X)$. Let $\cF = 
\{\, h \in \se L k\, \vert\,$ There is $g \in \cG$ such that $g \spez\, h\, \}$ be
the AOS-fan consisting of the depth-\y k successors of elements of\, \cG\, $(cf.\, \ref{fansS})$.
\und{Assume}:
\vspace{-0.15cm}

\nhp $(*)\;\;\;\; \fa\, h,h' \in \cF,\;\; {\mathrm {card}}\,(\{ g \in\, \cG\, \vert\, g
\spez\, h \}) = {\mathrm {card}}\,(\{ g \in\, \cG\, \vert\, g \spez\, h' \})$. \hfl
\vspace{-0.15cm}

\nhp Let $\cB = \{ \se h 1, \dots, \se h r \}$ be a basis of \cF\, $($as an
AOS\,$)$, and let\, \cC\ be a basis of the AOS-fan \lbr
$ \se P {h_1} = \{ g \in\, \cG\, \vert\, g \spez\, \se h 1 \}$ $($see \ref{permsetsS+C}\,(1). 
For\, $i = 2, \dots, r$, let $\se g i \in\, \cG$ be such that  $\se g i \spez\, \se h i$.
\vspace{-0.15cm}

\noi \und{Then},\; $\cC\, \cup \{ \se g 2, \dots, \se g r \}$ is a basis of\, \cG. 
\epr
\vspace{-0.35cm}

\nhp {\bf Proof.} If \y r = 1, then $\cF = \cB = \{\se h 1\}$, whence $\cG = \{g
\in \cG\, \vert\, g \spez\, \se h 1 \}$, and the result holds by the choice of
\cC. Henceforth we assume $r \geq 2$. We observe:
\vspace{-0.15cm}

\nhp ---\; $r = {\mathrm {card}}(\cB) = {\mathrm {dim}}(\cF)$. Since \cF\ is an AOS-fan,\, 
${\mathrm {card}}(\cF) = \po {r - 1} 2$.
\vspace{-0.15cm}

\nhp ---\; For every $h \in \cF$, $\se A h = \{g \in \cG\, \vert\, g \spez\, h
\}$ is a AOS-fan; this follows from the assumption that \cG\, is an AOS-fan,
since \se A h is closed under the product of any three of its elements, cf. Lemma
\ref{moreprods}\,(b).
\vspace{-0.15cm}

\nhp ---\; $\se A h \cap\, \se A {h'} = \0$\, for $h \neq h'$\, in\, \cF.
\vspace{-0.15cm}

\nhp By assumption $(*)$, card\,(\se A h) = card\,(\se A {h'}) (= $\po {p -
1} 2$, say), for $h, h' \in \cF$. Since $\cG = \bigcup_{h \in \cal F} \se A h$, we
get card\,(\cG) = card\,(\cF)$\cdot$\,card\,(\se A h) (any $h \in \cF$), and
then card\,(\cG) = $\po {r-1} 2 \cdot\, \po {p-1} 2 = \po {p+r-2} 2$; hence
dim\,(\cG) = $p+r-1$. Since card\,($\,\cC\, \cup \{ \se g 2, \dots, \se g r \}$)
= $p+r-1$, it suffices to prove:
\vspace{-0.15cm}

\nhp $(**)\;\;\; \cC\, \cup \{ \se g 2, \dots, \se g r \}$\, is an independent
set.
\vspace{-0.15cm}

\nhp \und{Proof of $(**)$}. Assume false.
\vspace{-0.15cm}

\nhp \und{Case 1.}\; Some $\se g {i_0}$, with $2 \leq \se i 0 \leq r$, is
dependent on the rest, i.e., there are\, $\cC'\, \sub\; \cC$\, and\, $J\, \sub\,
\{2, \dots, r\} \setminus \{\se i 0\}$ so that\, $\se g {i_0} = \prod_{c\in
\cal C'} c \cdot\, \prod_{j\in J} \se g j$,\, i.e.,
\vspace{-0.2cm}

\nhp $(+)\;\;\; \prod_{c\in \cal C'} c = \prod_{j\in J\cup\{i_0\}} \se g j.$
\vspace{-0.2cm}

\nhp --- If card\,($\cC'$) is odd, since \se A {h_1} is an AOS-fan, and hence a
closed set, the left-hand side of $(+)$ is an element $g' \spez\, \se h 1$,
and we have\, $g' \cdot\, \prod_{j\in J\cup\{i_0\}} \se g j = 1$. Setting $A =
\{g'\}\, \cup\, \{ \se g j\, \vert\,$ $j\in J\cup\{i_0\}\,\}\, \sub\, \se L
{k+1}$ and $B = \{\se h 1\}\, \cup\, \{ \se h j\, \vert\, j\in
J\cup\{i_0\}\,\}\, \sub\, \se L k$, the assumptions of Lemma \ref{dep-fans} are
met. Since \y A is dependent, so is \y B, contradicting that $B\, \sub\, \cB$
and \cB\, is a basis of \cF, whence an independent set.
\vspace{-0.15cm}

\nhp --- If $\cC'$ = \0, the same argument works, yielding a contradiction.
\vspace{-0.15cm}

\nhp --- Assume card\,($\cC'$) even $>\, 0$. Fix $\se c 0 \in \cC'$. Then
card\,($\cC' \setminus \{\se c 0\}$) = odd, and $g' = \prod_{c\in {\cal C'} \setminus
\{ c_0 \}} c \in \se L {k+1}$; also $g' \spez\, \se h 1$, and we have:
\vspace{-0.3cm}

\nhp \hfl $\se c 0 \cdot\, g' \cdot\, \prod_{j\in J\cup\{i_0\}} \se g j = 1.$ \hfl
\vspace{-0.1cm}

\nhp Pick any index $\se j 0 \in J\,\cup\, \{i_0\}$ (so, $\se j 0 \geq 2$).
Since $\se c 0\,,\, g' \spez\, \se h 1$ and $\se g {j_0} \spez\, \se h {j_0}$, 
Lemma \ref{moreprods}(b) yields $\se {g'} {j_0} := \se c 0\, g'\,
\se g {j_0} \spez\, \se h {j_0}$, and $\se {g'} {j_0} \cdot\,
\prod_{j\in (J\cup\{i_0\})\setminus \{j_0\}} \se g j = 1$. Hence, 
$A = \{\se {g'} {j_0}\}\, \cup\, \{\se g j\, \vert$ 
$j \in (J\, \cup\, \{\se i 0\})\setminus \{\se j 0\} \}$ is a dependent subset
of \se L {k+1}. Setting $B = \{\se h j\, \vert\, j \in J\, \cup\, \{\se i 0\}\}$
the assumptions of Lemma \ref{dep-fans} are met, and hence \y B is also
dependent, contradicting that $B\, \sub\, \cB$.
\vspace{-0.15cm}

\nhp \und{Case 2.}\; Some $\se c 0 \in \cC$ is dependent on the rest.
\vspace{-0.15cm}

\nhp Then, there are\, $\cC'\; \sub\; \cC \setminus \{\se c 0\}$\, and\, $J\,
\sub\, \{2, \dots, r\}$ so that
\vspace{-0.2cm}

\nhp $(++)\;\;\; \se c 0 = \prod_{c\in \cal C'} c \cdot \prod_{j\in J} \se g j.$
\vspace{-0.1cm}

\nhp Note that $J \neq \0$ (otherwise\, \cC\, would be dependent). Taking \y J
minimal so that $(++)$ holds, and picking $\se j 0 \in J$, it follows that\, \se
c 0\, is not in the closure of\, $\cC'\, \cup\, \{\se g j\, \vert$ $j \in J \setminus \{\se j 0\} \}$
(cf. Warning, end of \ref{notatff}\,(b)). By the exchange property,\, 
$\se g {j_0}$ is in the closure of\, $\cC'\, \cup\, \{\se c 0\}\, \cup\, \{\se g
j\, \vert\, j \in J \setminus \{\se j 0\} \}$, contrary to the result of Case 1.
\hfl $\Box$

\vspace{-0.35cm}

\bct \label{stgen} {\bf Standard generating systems.}
\vspace{-0.1cm}

For any finite ARS-fan, $(X,F)$, we will construct, by induction on \y k, $1
\leq k \leq n = \ell(X)$, a class of bases $\se {\cB} k$ of the AOS-fan $\se L
k(X)$. Each basis $\se {\cB} k$ will be required to satisfy the additional
condition:
\vspace{-0.2cm}

\nhp $(*) \;\;$ For $k \leq j \leq n$,\, $\se {\cB} k \, \cap\, \se {S^k} j$ is
a basis of the AOS-fan $\se {S^k} j$. 
\vspace{-0.2cm}

\nhp This additional requirement guarantees that the inductive construction of
the $\se {\cB} k$'s is not interrupted before the \y n-th (and last) step.
The construction uses Proposition \ref{choicebas} and the results from \S\,\ref{involut} 
above. The set $\cB = \bigcup_{k=1}^n \se {\cB} k$ is called a 
{\bf standard generating system} for $(X,F)$.
\vspace{-0.1cm}

\noi {\bf Construction of standard generating systems.}
\vspace{-0.1cm}

\nhp \und{Level 1}. It suffices to observe that a basis $\se {\cB} 1$ of \se L 1
satisfying condition $(*)$ exists. Begin by choosing a basis $\se {\cB} 1(n)$ of
the AOS-fan $\se {S^1} n = \se {C^1} n$ (cf. Corollary \ref{fansS}). $\se {S^1} n$ 
is a closed subset (cf. Warning, end of \ref{notatff}\,(b)) of the (AOS-)fan 
$\se {S^1} {n-1} = \se {S^1} n  \, \cup\, \se {C^1} {n-1}$; hence,\, $\se {\cB} 1(n)$
is an independent subset of $\se {S^1} {n-1}$; choose $\se {\cB} 1(n-1)$ to be a
basis of $\se {S^1} {n-1}$ extending $\se {\cB} 1(n)$.
\vspace{-0.15cm}

In general, assume that, for $1 < j \leq n$ an increasing sequence $\se {\cB}
1(n)\, \sub \dots \sub\, \se {\cB} 1(j)$ of independent subsets of \se L 1 has
been chosen so that $\se {\cB} 1(\ell)$ is a basis of the AOS-fan\; $\se {S^1}
\ell\;\, (j \leq \ell \leq n)$. As above, $\se {\cB} 1(j)$ is an independent
subset of the fan $\se {S^1} {j-1} = \se {S^1} j  \, \cup\, \se {C^1} {j-1}$.
Let $\se {\cB} 1(j-1)$ be a basis of $\se {S^1} {j-1}$ extending $\se {\cB}
1(j)$. Set $\se {\cB} 1 = \se {\cB} 1(1)$; by construction, $\se {\cB} 1 \,
\cap\, \se {S^1} j = \se {\cB} 1(j)$ is a basis of $\se {S^1} j$.
\vspace{-0.15cm}

\nhp \und{Induction step.} Given an integer \y k, $1 \leq k < n$, assume there
exists a basis\, $\se {\cB} k$ of \se L k satisfying property $(*)$; thus, for
$k \leq j \leq n$, $\se {\cB} k(j) = \se {\cB} k \cap\, \se {S^k} j$ is a basis
of $\se {S^k} j$. Further, since $\se {S^k} n\, \sub \dots \sub\, \se {S^k} k =
\se L k$, we have $\se {\cB} k(n)\, \sub \dots \sub\, \se {\cB} k(k) = \se {\cB}
k$. Using Proposition \ref{choicebas} we define a\vspace{0.08cm} subset 
$\se {\cB} {k+1}$ of $\se L {k+1}$ as follows:

\nhp --- Firstly, fix an element $\se h 0 \in \se {\cB} k(n)$ (this set is
non-empty because $n = \ell(X)$). Pick a basis $\se {\cB} {k+1}(n, \se h 0)$ of
the AOS-fan $\{g \in \se{S^{k+1}} n\, \vert\; g \spez\, \se h 0 \}$.
\vspace{-0.15cm}

\nhp --- Next, for each $h \in (\se {\cB} k \, \cap\, \se {S^k} {k+1}) \setminus
\{\se h 0\}$ there is a maximal index \y j $ = j(h)$,\; $k + 1 \leq j \leq n$,\, so that
$h \in \se {\cB} k \, \cap\, \se {S^k} j = \se {\cB} k(j)$; clearly, $h\ \nen\
\se {S^k} {j+1}$, whence $h \in \se {C^k} j = \se {S^k} j \setminus \se {S^k}
{j+1}$ (if \y j = \y n, then $h \in \se {S^k} n = \se {C^k} n$). Since $j \geq k
+ 1$, we have $\{g \in \se{C^{k+1}} j\, \vert\; g \spez\, h \} \neq \0$.
Choose an element $\se g h \in \se{C^{k+1}} j$ so that $\se g h \spez\, h$.
\vspace{-0.15cm}

\nhp --- Finally, set 
\vspace{-0.2cm}

\noi [*]\;\; $\se {\cB} {k+1} = \se {\cB} {k+1}(n, \se h 0)\, \cup\,
\{\se g h\, \vert\, h \in (\se {\cB} k \, \cap\, \se {S^k} {k+1}) \setminus \{\se h 0\} \}$.
\vspace{-0.1cm}

\noi {\bf Claim.} For $k + 1 \leq p \leq n,\; \se {\cB} {k+1}\, \cap\, \se{S^{k+1}} p$ 
is a basis of $\se{S^{k+1}} p$.
\vspace{-0.15cm}

\noi {\bf Proof of Claim.} We apply Proposition \ref{choicebas} with the following choice
of parameters:
\vspace{-0.15cm}

\noi --- $\cG = \se{S^{k+1}} p$ (whence $\cF = \se {S^k} p$, since $k + 1 \leq p$);
\vspace{-0.25cm}

\noi --- $\cB = \se {\cB} k \, \cap\, \se {S^k} p$ (a basis of\, \cF);
\vspace{-0.25cm}

\noi --- $\cC = \se {\cB} {k+1}(n, \se h 0)$ (a basis of\, $\{ g \in \se{S^{k+1}} n \sth g \spez \se h 0)$).
\vspace{-0.1cm}

\noi Proposition \ref{cardsetsS+C}\,(a) shows that the cardinality assumption
\vspace{-0.1cm}

\nhp \hfl card\,($\{g \in \se{S^{k+1}} p\, \vert\; g \spez\, h \}$) =
card\,($\{g \in \se{S^{k+1}} p\, \vert\; g \spez\, h' \}$),\;\;\; $(h, h' \in
\se {S^k} j)$ \hfl
\vspace{-0.15cm}

\nhp of \ref{choicebas} holds. We conclude that 
\vspace{-0.15cm}

\noi \hfl $\cD := \se {\cB} {k+1}(n, \se h 0)\, \cup\, \{\se g h\, \vert\, h \in
(\se {\cB} k \, \cap\, \se {S^k} p) \setminus \{\se h 0\} \}$ \hfl
\vspace{-0.2cm}

\noi is a basis of $\se{S^{k+1}} j$. The Claim follows from:
\vspace{-0.2cm}

\noi (\dag)\;\; $\se {\cB} {k+1}\, \cap\, \se{S^{k+1}} p = \;\cD$.
\vspace{-0.1cm}

\noi \und{\it Proof of $(\dag)$}. Since $\se {\cB} {k+1}(n, \se h 0)\, \sub\, \cD\, \cap \se {\cB} {k+1}$
(see [*]), we need only prove:
\vspace{-0.1cm}

\noi ($\sub$)\; If $h \in (\se {\cB} k \, \cap\, \se {S^k} {k+1}) \setminus \{\se h 0\} \}$\, and\,
$\se g h \in \se{S^{k+1}} p$, then\, $h \in \se {\cB} k \, \cap\, \se {S^k} p$.
\vspace{-0.2cm}

\noi This clearly follows from $\se g h \in \se{S^{k+1}} p\,,\, \se g h \spez h$\, and\, 
$h \in \se {S^k} {k+1}$.
\vspace{-0.15cm}

\noi ($\supseteq$)\; Since $k + 1 \leq p$, we have\, 
$\se {\cB} {k}\, \cap\, \se{S^{k}} p = \se {\cB} {k}(p)\, \sub\, \se {\cB} {k}(k+1) =
\se {\cB} k \, \cap\, \se {S^k} {k+1}$. On the other hand, if 
$h \in (\se {\cB} k \, \cap\, \se {S^k} p) \setminus \{\se h 0\}$
and, as above, $j(h)$ denotes the largest index $j$ so that $k + 1 \leq j \leq n$ and 
$h \in \se {\cB} k(j)$, we have $p \leq j(h)$, whence $\se{S^{k+1}} {j(h)}\, \sub\, \se{S^{k+1}} p$.
By choice, $\se g h \in \se{C^{k+1}} {j(h)}$; it follows that $\se g h \in \se{S^{k+1}} p$, as
required. \hfl $\Box$
\ect
\vspace{-0.65cm}

\bres (a) In general, there are many different standard generating systems 
for a finite ARS-fan $(X,F)$. The construction in \ref{stgen} allows for 
several choices of the bases $\se {\cB} 1(j)$ $(1 \leq j \leq n)$ and, at 
each successive step, \y k, for many choices of elements 
$\se h 0 \in \se {\cB} k(n)$, of bases $\se {\cB} {k+1}(n, \se h 0)$, and 
of elements $\se g h \in \se{C^{k+1}} n$ under each $h \in (\se {\cB} k \, 
\cap\, \se {S^k} {k+1}) \setminus \{\se h 0\}$. In spite of this lack of 
uniqueness, we shall prove below that any standard generating system 
determines the isomorphism type of a finite ARS-fan.
\vspace{-0.2cm}

\nhp (b) Some of the\vspace{-0.05cm} sets $\se {C^k} j = \se {S^k} j \setminus \se {S^k} {j+1}$ 
may be empty. However, if $\se {C^k} j \neq \0$, then, necessarily,\, 
$\se {\cB} k \, \cap\, \se {C^k} j \neq \0$. Indeed, if \y j = \y n, then 
$\se {C^k} n \neq \0$ (as $n = \ell(X)$) and $\se {C^k} n = \se {S^k} n$ 
is an AOS-fan; since $\se {\cB} k \, \cap\, \se {C^k} n$ is a basis of 
$\se {C^k} n$, it must contain at least one element. If $j < n$, since 
$\se {S^k} {j+1}$ is a fan, it is a closed set; as it is disjoint from 
$\se {C^k} j$, then no\vspace*{-0.05cm} element of $\se {C^k} j$ is dependent on 
$\se {S^k} {j+1}$. Hence, any basis of $\se {S^k} j = \se {S^k} {j+1}\, \cup\, 
\se {C^k} j$ must contain an element of $\se {C^k} j$. \hfl $\Box$
\eres

\vspace{-0.45cm}

Any standard generating system for a finite ARS-fan has the following property:

\vspace{-0.4cm}

\bco \label{stgenprop} Let\, \cB\, be a standard generating system for a finite
ARS-fan $(X,F)$; let $n = \ell(X)$, and $1 \leq k \leq m \leq n$. Then, for
every\, $g \in \se {\cB} m = \cB\, \cap\, \se L m(X)$, the\vspace{0.05cm} unique depth-\y k
successor of\, \y g in \y X belongs to\, \cB\, $($hence to $\se {\cB} k = \cB\,
\cap\, \se L k(X))$.
\eco

\vspace{-0.35cm}
\nhp {\bf Proof.} By the construction in \ref{stgen} this holds for $m = k + 1$. In fact,
let $g \in \se {\cB} {k+1}\,,\, h' \in \se L k$ and $g \spez h'$. By the definition of
$\se {\cB} {k+1}$ and uniqueness of the successor of $g$ in $\se L k$, if 
$g \in \se {\cB} {k+1}(n, \se h 0)$, then $h' = \se h 0 \in \se {\cB} k(n)\, \sub\, \se {\cB} k$;
if $g = \se g h$, with $h \in (\se {\cB} k \, \cap\, \se {S^k} {k+1}) \setminus \{\se h 0\}$,
we get $h' = h \in \se {\cB} k$. Then iterate. \hfl $\Box$

For the proof of the Isomorphism Theorem \ref{isothffans} below we shall need the
characterizations of ARS-morphisms between fans proved in \ref{fan-morph} and 
\ref{fan-maps} below, which, in turn, follow from the Small Representation Theorem 
\ref{smallrep}.

\vspace{-0.3cm}

\bdf \label{pres3-prod} Let $G,H$ be RSs, let $\se X G, \se X H$ be their 
character spaces, and let $Z\, \sub\, \se X G$\,. A map $F: Z\, \lra\, \se X H$  
{\bf preserves 3-products} $($in \y Z$)$ iff for all $\se h 1, \se h 2, \se h 3 \in Z$,

\nhp \hfl  $\se h 1 \se h 2 \se h 3 \in Z\;\; \Ra\;\; F(\se h 1 \se h 2 \se h 3) = 
F(\se h 1) F(\se h 2) F(\se h 3).$ \hfl $\Box$
\edf

\vspace{-0.8cm}

\bpr \label{smallrep} {\em (Small representation theorem)}. Let $G$ be a RS. 
The following conditions are equivalent for a map\, $f: \se X G\, \lra\, \bf 3$\,:
\vspace{-0.15cm}

\nhp $(1)$ a$)$ \y f is continuous in the constructible topology of $\se X G$.
\vspace{-0.15cm}

\nhp \;\;\;\;\, b$)$ \y f preserves $3$-products in $\se X G$.
\vspace{-0.15cm}

\nhp $(2)$ \y f is represented by an element of \y G: there is $a \in G$ so that $f = \h a$.
\epr

\vspace{-0.45cm} 

\nhp [\,$\h a: \se X G \lra\ \bf 3$\, denotes ``evaluation at \y a": 
for $h \in \se X G$,\; $\h a(h) = h(a)$.]

\nhp {\bf Proof.} (2)\,$\Ra$\,(1) is clear since the evaluation maps have 
properties (1.a) and (1.b).
\vspace{-0.15cm}

\nhp (1)\,$\Ra$\,(2). We use the representation theorem \cite{M}, Cor. 8.3.6, 
p. 162. It suffices to check that the assumptions of this theorem as well as 
one of the equivalent conditions in its conclusion hold under our hypotheses in 
(1). In our notation, the conditions to be checked are: for $x, y \in \se X G$,

\vspace{-0.15cm}

\nhp ($\dag$)\;\;\;\, $f(x) = 0$\, and\, $Z(x)\, \sub\, Z(y)$\; implies\; $f(y) = 0$. 
\vspace{-0.15cm}

\nhp ($\dag\dag$)\;\; $f(x) \neq 0$\, and\, $x^{-1}[0,1] \supseteq y^{-1}[0,1]$\; 
implies\; $f(x) = f(y)$.
\vspace{-0.15cm}

\nhp ($\dag\dag\dag$)\, For any saturated prime ideal \y I of \y G, either 

\vspace{-0.2cm}

\nhp \hspace{0.9cm} (i)\;\, $f \lceil \{u \in \se X G\, \vert\, Z(u) = I\} = 0$,\, or 

\vspace{-0.15cm}

\nhp \hspace{0.9cm} (ii)\; $\prod_{i=1}^{4}f(\se x i) = 1$\, for any 4-element 
AOS-fan $\{\se x 1,\dots \se x 4\}$\, in\, $\{u \in \se X G\, \vert\, Z(u) = I\}$. 

\nhp --- Condition ($\dag$) follows at once from Lemma \ref{char-zeroset}\,(2) 
(as $Z(x)\, \sub\, Z(y)\;\Ra\; y = yx^2$).
\vspace{-0.15cm}

\nhp --- Condition ($\dag\dag$) follows from Lemma \ref{char-specializ}\,(3),(5):\; 
$x^{-1}[0,1] \supseteq y^{-1}[0,1]\; \Ra\; x = x^2y$. Since $f(x) \neq 0\; \Ra\; 
f(x^2) = 1$, assumption (1.b) implies\; $f(x) = f(x^2)f(y) = f(y)$.
\vspace{-0.15cm}

\nhp --- As for ($\dag\dag\dag$), if (i) does not hold, ($\dag$) implies\, 
$f(u) \neq\, 0$\, for all\, $u \in \se X G$\, such that\, $Z(u) = I$. Let\, 
$\{\se x 1,\dots \se x 4\}$\, be an AOS-fan in $\{u \in \se X G\, \vert\, Z(u) = I\}$. 
Thus,\, $\se x 4 = \se x 1\se x 2\se x 3$\, and\, $f(\se x i) \neq 0$ for 
$i = 1, \dots, 4$. Assumption (1.b) gives $f(\se x 4) = f(\se x 1)f(\se x 2)f(\se x 3) 
\neq\, 0$, i.e.,\, $\prod_{i=1}^{4}f(\se x i) = 1$. \hfl $\Box$ 
\vspace{-0.3cm}

\bco \label{fan-morph} A map $F: (\se X 1,\se F 1)\, \lra\, (\se X 2,\se F 2)$
between ARS-fans is an ARS-morphism iff\, \y F is continuous for the
constructible topology $($of both source and target\,$)$ and preserves
$3$-products in \se X 1 $($cf. $\ref{pres3-prod}$\,$)$.
\eco
\vspace{-0.35cm}

\nhp {\bf Proof.} ($\La$) If \y F has the stated properties and $a \in \se F 2$,
then $\h a\; \com\, F: \se X 1\, \lra\, \bf 3$ also has those properties, and,
by Proposition \ref{smallrep}, is represented by an element of \se F 1; hence,\,
\y F is an ARS-morphism (cf. \ref{levels}\,(c.i)). 
\vspace{-0.15cm}

\nhp ($\Ra$) Assume\, \y F is an ARS-morphism. For continuity it suffices to
show that $F^{-1}[V]$ is open constructible in \se X 1 whenever \y V is basic
open for the constructible topology of \se X 2, i.e., of the form $V = U(\se a
1, \dots, \se a n)\, \cap\, Z(a)$ with $a, \se a 1, \dots, \se a n \in \se F 2$
(see \cite{M}, p. 111).
By the assumption on\, \y F, there are\; $b, \se b 1, \dots, \se b n \in \se F
1$\; such that\, $\h a\; \com\, F = \h b$\, and\, $\h {\se a i}\; \com\, F = \h
{\se b i}$\, for\, $i = 1, \dots, n$. These functional identities imply
$F^{-1}[V] = U(\se b 1, \dots, \se b n)\, \cap\, Z(b)$, as required.

\vspace{-0.15cm}

Preservation of 3-products by\, \y F\, follows\vspace{0.05cm} easily from the same property
for\, $\h a$\, and\, $\h b$\, using the functional identity\, $\h a\; \com\, F =
\h b$. \hfl $\Box$

\vspace{-0.35cm}

\ble \label{fan-maps} Let $(\se X 1,\se F 1),\, (\se X 2,\se F 2)$ be ARS-fans.
\vspace{-0.15cm}

\nhp $(1)$ For a map $F: \se X 1\, \lra\, \se X 2$ the following are equivalent:
\vspace{-0.15cm}

\tab{\nhp $(1)$} $($i\,$)$\; \y F preserves $3$-products in \se X 1.
\vspace{-0.15cm}

\tab{\nhp $(1)$} $($ii\,$)$ a$)$ \y F preserves $3$-products of elements {\em of the same
level}: for all $I \in$ {\em Spec}$(\se F 1)$\, and all\vspace{0.08cm} \hspace*{1.8cm} $\se h 1, 
\se h 2, \se h 3 \in \se L I(\se X 1)$,\; $F(\se h 1 \se h 2 \se h 3) = 
F(\se h 1)F(\se h 2)F(\se h 3).$ 
\vspace{-0.15cm}

\tab{\nhp $(1)$ $($ii\,$)$\;} b$)$ \y F is monotone for the specialization order: for\, 
$g,h \in \se X 1$,\, $g\,\spezz 1 \, h \;\; \Ra \;\; F(g)\, \spezz 2 \, F(h)$.

\vspace{-0.2cm}

\nhp $(\spezz i$ denotes specialization in \se X i$)$.
\vspace{-0.15cm}

\nhp $(2)$ If $(\se X 1,\se F 1)$ is finite, any map verifying one of the
equivalent conditions $($i$)$ or $($ii$)$ in $(1)$ is a morphism of ARSs.
\vspace{-0.15cm}

\nhp $(3)$ If both $(\se X 1,\se F 1),\, (\se X 2,\se F 2)$ are finite, any
bijection $F: \se X 1\, \lra\, \se X 2$ verifiying one of the equivalent
conditions in $(1)$ is an isomorphism of ARSs.
\ele

\vspace{-0.35cm}

\nhp {\bf Proof.} (1). (i) $\Ra$ (ii). (ii.a) is a special case of (i).
\vspace{-0.15cm}

\nhp (ii.b)\;\; $g\,\spezz 1\, h \;\; \Lra \;\; h = h^2g$ (Lemma
\ref{char-specializ}). By (i), $F(h) = F(h)^2F(g)$, and\vspace*{0.05cm} this equality (in \se X 2) 
gives $F(g)\, \spezz 2\, F(h)$.
\vspace{-0.15cm}

\nhp (ii) $\Ra$ (i). Let $\se h 1, \se h 2, \se h 3$ be any three elements in
\se X 1; say 
$Z(\se h 3)\, \sub\, Z(\se h 2)\, \sub\, Z(\se h 1)$. Let $I = Z(\se h 1)$ and
for\, $i = 2, 3$\, let\, $\se {h'} i$\, be the unique successor of\, \se h i\,
in\, $\se L I(\se X 1)$; Lemma \ref{prods} shows that\, $\se h 1 \se h 2 \se h 3
= \se h 1 \se {h'} 2 \se {h'} 3$; then, assumption (ii.a) gives
\vspace{-0.15cm}

\nhp \hfl $F(\se h 1 \se h 2 \se h 3) = F(\se h 1) F(\se {h'} 2) F(\se {h'} 3).$ \hfl
\vspace{-0.15cm}

\nhp By (ii.b) we have\, $F(\se h i)\, \spezz 2\, F(h'_i),\; (i = 2,3)$. 
Next, note that $Z(F(h'_i))\, \sub\, Z(F(\se h 1))$ for $i = 2,3$. In fact, 
since $\se h 1, \se {h'} 2$ belong to the same level $\se L I$, $Z(\se h 1) = Z(\se {h'} 2)$, 
and Lemma \ref{char-zeroset}\,(2) yields $\se {h^2} 1 = {\se {h'} 2}^2$; scaling by    
$\se h 1$ gives $\se h 1 = \se h 1 {\se {h'} 2}^2$. Since $F$ preserves 3-products of 
the same level, $F(\se h 1) = F(\se h 1) F(\se {h'} 2)^2$ which, by \ref{char-zeroset}\,(1),
gives $Z(F(\se {h'} 2))\, \sub\, Z(F(\se h 1))$. Same argument for $i = 3$. 

\vspace{-0.15cm}

\noi Using \ref{prods} again, these inclusions and $F(\se h i)\, \spezz 2\, F(\se {h'} i),\; (i = 2,3)$ 
prove:
\vspace{-0.15cm}

\nhp \hfl $F(\se h 1) F(\se {h'} 2) F(\se {h'} 3) = F(\se h 1) F(\se h 2) F(\se
h 3),$ \hfl
\vspace{-0.2cm}

\nhp as required.

\nhp (2) follows at once from Corollary \ref{fan-morph}, since the continuity
requirement is automatically fulfilled in this case: the constructible topology
in\, \se X 1 is discrete.
\vspace{-0.15cm}

\nhp (3) By (2) it only remains to prove that\, $F^{-1}: \se X 2\, \lra\, \se X
1$\, preserves 3-products in\, \se X 2. Let\, $\se g 1, \se g 2, \se g 3 \in \se
X 2$\, and let $\se h i = F^{-1}(\se g i), \; i = 1, 2, 3$. From (1.i) we have\,
$F(\se h 1 \se h 2 \se h 3) = \se g 1 \se g 2 \se g 3$. Composing both sides of
this equality with\, $F^{-1}$\, gives the desired conclusion:
\vspace{-0.15cm}

\nhp \hfl $F^{-1}(\se g 1 \se g 2 \se g 3) = F^{-1}(F(\se h 1 \se h 2 \se h 3))
= \se h 1 \se h 2 \se h 3 = F^{-1}(\se g 1) F^{-1}(\se g 2) F^{-1}(\se g 3).$
\hfl \;\;\; $\Box$
  
\nhp {\bf Remark.} Note that any isomorphism of ARS-fans preserves depth.
\vspace{-0.3cm}

\bth \label{isothffans} {\em (The isomorphism theorem for finite ARS-fans.)} Let
$(\se X 1, \se F 1),\, (\se X 2,\se F 2)$ be finite ARS-fans and let
$\spezz 1, \spezz 2$ denote their respective specialization orders. If\, 
$(\se X 1, \spezz 1)$ and $(\se X 2, \spezz 2)$ are order-isomorphic, then 
$\se X 1$ and $\se X 2$ are isomorphic ARSs.
\eth

\vspace{-0.35cm}
\nhp {\bf Proof.} The order-isomorphism assumption implies:
\vspace{-0.15cm}

\nhp (1) $\ell(\se X 1) = \ell(\se X 2)$ ( = \y n, say, fixed throughout the
proof).
\vspace{-0.15cm}

\nhp (2) For $1 \leq k \leq j \leq n$, card($\se {C^k} j(\se X 1))$ = card($\se
{C^k} j(\se X 2))$.
\vspace{-0.15cm}

\nhp The proof of (2) is an easy exercise. Since $\se {C^k} {\ell} \cap \se
{C^k} {\ell'} = \0$ for $k \leq \ell \neq \ell' \leq n$ and $\se {S^k} j =
\bigcup_{\ell=j}^{\,n} \se {C^k} {\ell}$, we get:
\vspace{-0.15cm}

\nhp (3) For $1 \leq k \leq j \leq n$,\; ${\mathrm {card}}(\se {S^k} j(\se X 1)) =\,
{\mathrm {card}}(\se {S^k} j(\se X 2))$.
\vspace{-0.15cm}

\nhp (4) For $1 \leq k < n$ and all $h \in \se {S^k} n(\se X 1)$, $h' \in \se
{S^k} n(\se X 2)$, we have:
\vspace{-0.15cm}

\nhp \hfl ${\mathrm {card}}(\{g \in \se {S^{k+1}} n(\se X 1)\, \vert\; g\, \spezz 1\,
h \}) =\, {\mathrm {card}}(\{g'\, \in \se {S^{k+1}} n(\se X 2)\, \vert\; g' \spezz 2\, 
h' \}).$ \hfl
\vspace{-0.15cm}

\nhp \und{Proof of (4)}. Consider the two-variable formula in the language 
$\{\leq \}$ of order:
\vspace{-0.2cm}

\nhp \hfl $\fh(x,y)\; :=\; x \in \se {S^{k+1}} n\, \w\; x \leq y$. \hfl
\vspace{-0.2cm}

\nhp (It is left as an exercise for the reader to write a first-order formula in
$\{ \leq \}$ expressing the notion $x \in \se {S^{k+1}} n$; cf. \ref{setsS+C}.)
\vspace{-0.2cm}

\nhp If \si\ denotes the order isomorphism between $(\se X 1, \spezz 1)$ and 
$(\se X 2, \spezz 2)$, for $g,h \in \se X 1$ we have:
\vspace{-0.1cm}

\nhp \hfl $(\se X 1, \spezz 1) \models\, \fh[g,h] \;\; \Lra \;\; (\se
X 2, \spezz 2) \models\, \fh[\si(g),\si(h)].$ \hfl
\vspace{-0.15cm}

\nhp It follows that\, \si\, maps\, $\{g \in \se {S^{k+1}} n(\se X 1)\, \vert\;
g\, \spezz 1\, h \}$\, bijectively onto\,  $\{g' \in \se {S^{k+1}} n(\se
X 2)\, \vert\; g'\, \spezz 2\, \si(h) \}$. Now, if\, $h \in \se {S^k}
n(\se X 1)$, then\, $\si(h) \in \se {S^k} n(\se X 2)$. If\, $h' \in \se {S^k}
n(\se X 2)$, apply Proposition \ref{cardsetsS+C} with $\se h 1 = h'$ and $\se h
2 = \si(h)$ to conclude. \hfl $\Box$

Since the sets in item (4) are AOS-fans (Corollary \ref{fansS}), they have the
same dimension, i.e., any bases of each of them have the same cardinality. If\,
${\cal B}^1, {\cal B}^2$ are standard generating systems for \se X 1, \se X 2,
respectively, then ${\cal B}^i\, \cap\, \se {S^k} j(\se X i)$ is a basis of the
fan $\se {S^k} j(\se X i)$, for $1 \leq k \leq j \leq n$ and $i = 1,2$; from (3)
we get:
\vspace{-0.15cm}

\nhp (5) For $1 \leq k \leq j \leq n$,\; ${\mathrm {card\,}}({\cal B}^1 \cap\, \se {S^k} j(\se
X 1)) =\, {\mathrm {card\,}}({\cal B}^2 \cap\, \se {S^k} j(\se X 2))$.
\vspace{-0.2cm}

\nhp In particular, for $\se {S^k} k = \se L k$ we obtain:
\vspace{-0.2cm}

\nhp (6) If $1 \leq k \leq n$, then\, ${\mathrm {card\,}}(\se {{\cal B}^1} k) =\, 
{\mathrm {card\,}}(\se {{\cal B}^2} k)$\;\; (where $\se {{\cal B}^i} k = {\cal B}^i \cap\, 
\se L k(\se X i)$).
\vspace{-0.15cm}

Next, we choose an arbitrary standard generating system\, ${\cal B}^1$ for \se X
1. By induction on \y k, $1 \leq k \leq n$, we construct a
standard generating system\, ${\cal B}^2$ of \se X 2 (${\cal B}^2 =
\bigcup_{k=1}^{\,n} \se {{\cal B}^2} k$) and a map $\se f k: \se {{\cal B}^1}
k\, \lra\, \se {{\cal B}^2} k$ so that:
\vspace{-0.2cm}

\nhp (7) i) For $k \leq j \leq n$,\;\; $\se f k[\,{\cal B}^1 \cap\, \se {S^k} j(\se X 1)] 
=\, {\cal B}^2 \cap\, \se {S^k} j(\se X 2)$.
\vspace{-0.2cm}

\nhp \hspace{0.4cm} ii) If\, $k < n,\; g \in \se {{\cal B}^1} {k+1},\, h \in \se
{{\cal B}^1} k$\, and\, $g\, \spezz 1\, h$, then\, $\se f {k+1}(g)\, \spezz 2\, 
\se f k(h)$.
\vspace{-0.15cm}

\nhp \und{Construction of ${\cal B}^2$ and the maps \se f k}.
\vspace{-0.15cm}

\nhp \und{Level 1}. $\se {{\cal B}^2} 1$ is built as in the level 1 step in
\ref{stgen}; with notation therein, $\se f 1: \se {{\cal B}^1} 1\, \lra\, \se
{{\cal B}^2} 1$ is any bijection mapping\,
$\se {\cal B} 1(j)$\, onto\, $\se {\cal B} 2(j)$, for\, $1 \leq j \leq n$. Such
a bijection exists by (5) above $(k = 1)$.
\vspace{-0.2cm}

\nhp \und{Induction step}. Assume $\se {{\cal B}^2} 1, \dots, \se {{\cal B}^2}
k$ and $\se f 1, \dots, \se f k$ already constructed, so that:
\vspace{-0.2cm}

\nhp --- For $1 \leq j \leq k$ and $j \leq {\ell} \leq n$,\; $\se {{\cal B}^2} j
\cap\, \se {S^j} {\ell}(\se X 2)$\, is a basis of the AOS-fan $\se {S^j}
{\ell}(\se X 2)$\, and\, $\se f j[\,\se {{\cal B}^1} j \cap\, \se {S^j} {\ell}(\se
X 1)] =\, \se{{\cal B}^2} j\, \cap\, \se {S^j} {\ell}(\se X 2)$.
\vspace{-0.2cm}

\nhp --- Condition (7.ii) holds for all \y j\, such that $1 \leq j < k$.
\vspace{-0.1cm}

The basis $\se {{\cal B}^2} {k+1}$, and along with it the map\, $\se f {k+1}$,
are defined by performing the construction of the inductive step in \ref{stgen},
with the following choice of parameters:
\vspace{-0.15cm}

\nhp --- If\, $\se h 0 \in \se {{\cal B}^1} k\, \cap\, \se {S^k} n(\se X 1)$,
and\, $\se {{\cal B}^1} {k+1}(n, \se h 0)$\, is a basis of the (AOS-)fan 
$\{g \in \se {S^{k+1}} n(\se X 1)\, \vert\; g\, \spezz 1\, \se h 0 \}$, then take\,
$\se {{\cal B}^2} {k+1}(n, \se f k(\se h 0))$\, to be a basis of the fan\, 
$\{g' \in \se {S^{k+1}} n(\se X 2)\, \vert\; g'\, \spezz 2 \, \se f k(\se h 0)\}$. 
This is possible since $\se f k(\se h 0) \in \se {{\cal B}^2} k \cap\, \se
{S^k} n(\se X 2)$, by (7.i). Using item (4), we let 
$\se f {k+1} \lceil\, \se {{\cal B}^1} {k+1}(n, \se h 0)$ be a bijection between 
$\se {{\cal B}^1} {k+1}(n, \se h 0)$ and $\se {{\cal B}^2} {k+1}(n, \se f k(\se h 0))$.
\vspace{-0.15cm}

\nhp --- If\, $g \in \se {{\cal B}^1} {k+1} \cap\, \se {C^{k+1}} j(\se X 1)$\,
with $k+1 \leq j \leq n$, but\, $g \not\in \se {{\cal B}^1} {k+1}(n, \se h 0)$,
then, by the construction performed in the inductive step of \ref{stgen}, if \y
h\ is the unique depth-\y k successor of \y g, we have $h \in \se {{\cal B}^1}
k\, \cap\, \se {C^k} j(\se X 1),\, h \neq \se h 0$ and $g = \se g h$. In this
case choose any element\, $g'\, \spezz 2\, \se f k(h)$ such that $g' \in
\se {C^{k+1}} j(\se X 2)$, and set $\se f {k+1}(g) = g'$. This is possible
since\, $\se f k(h) \in \se {{\cal B}^2} k \cap\, \se {C^k} j(\se X 2)$ (which
follows easily from (7.i)). Clearly, this construction guarantees that (7.i) and 
(7.ii) hold for $k + 1$.
\vspace{-0.15cm}

Note that (7.ii) implies, by iteration, its own generalization:
\vspace{-0.1cm}

\nhp (7) iii) If\, $1 \leq k < m \leq n,\, g \in \se {{\cal B}^1} m,\, h \in \se
{{\cal B}^1} k$\, and\, $g\, \spezz 1\, h$, then\, $\se f m(g)\, \spezz 2\, \se f k(h)$.
\vspace{-0.15cm}

Since $\se {{\cal B}^i} k = {\cal B}^i \cap\, \se L k(\se X i)$ is a basis of
the AOS-fan $\se L k(\se X i)$, $i = 1,2$, we get:
\vspace{-0.2cm}

\nhp (8) The bijection\, \se f k\, extends (uniquely) to an AOS-isomorphism\,
$\td{\se f k}: \se L k(\se X 1)\, \lra\, \se L k(\se X 2)$\, mapping $\se
{S^k} j(\se X 1)$ onto $\se {S^k} j(\se X 2)$, for all \y j\ such that $k \leq j \leq n$.
\vspace{-0.2cm}

Now set\, $F: \se X 1\, \lra\, \se X 2$\, to be\, $F =\, \bigcup_{k=1}^{\,n}\td{\se f k}$. We prove:
\vspace{-0.1cm}

\nhp {\bf Claim.} \y F is an isomorphism of ARSs.
\vspace{-0.1cm}

\nhp {\bf Proof of Claim.}\, Since\, $\se X i =\, \bigcup_{k=1}^{\,n} \se L
k(\se X i)$ (disjoint union) for $i = 1,2$, and $\td{\se f k}$ maps $\se L
k(\se X 1)$ bijectively onto $\se L k(\se X 2)$, we have:
\vspace{-0.15cm}

\nhp (a) \y F\ is well-defined and bijective.
\vspace{-0.15cm}

\nhp (b) For all \y k, $1 \leq k \leq n$, \y F preserves 3-products \und{in \se L k}.
\vspace{-0.15cm}

\nhp This is clear: by (8)\, $F \lceil\, \se L k(\se X 1) = \td{\se f k} :
\se L k(\se X 1)\, \lra\, \se L k(\se X 2)$\, is an isomorphism of AOS-fans.
\vspace{-0.1cm}

\nhp (c) \y F is monotone for the specialization order.
\vspace{-0.15cm}

\nhp Let $g,h \in \se X 1$ be such that $g\, \spezz 1\, h$; say $d(g) = m
\geq d(h) = k$. We must prove $F(g)\, \spezz 2\, F(h)$. Since\, $\se
{{\cal B}^1} m$\, generates\, $\se L m(\se X 1)$, then\, $g = \se g 1 \cdot\
\dots\ \cdot \se g r$\ with\, $\se g 1, \dots, \se g r \in \se {{\cal B}^1} m$\,
and \y r necessarily odd (possibly = 1). By Corollary \ref{stgenprop}, if\, \se
h i is the unique depth-\y k successor of\, \se g i, then $\se h i \in \se
{{\cal B}^1} k$. Also, $\se g i\, \spezz 1\, \se h i$\; ($i = 1, \dots,
r$)\; implies\, $g = \se g 1 \cdot\ \dots\ \cdot \se g r\, \spezz 1\, \se
h 1 \cdot\ \dots\ \cdot \se h r$\, (\ref{moreprods}\,(a)).\vspace{0.05cm} Since both\, 
\y h\, and\, $\se h 1 \cdot\ \dots\ \cdot \se h r$\, are successors of \y g\, of the
same level \y k, we get\, $h = \se h 1 \cdot\ \dots\ \cdot \se h r$. As\, \y F\,
preserves products of any odd number of elements of the same level, we have:
\vspace{-0.15cm}

\nhp \hfl $F(g) = F(\se g 1) \cdot\ \dots\ \cdot F(\se g r)$ \;\;\; and \;\;\;
$F(h) = F(\se h 1) \cdot\ \dots\ \cdot F(\se h r)$. \hfl
\vspace{-0.2cm}

\nhp Since\, $\se g i\, \spezz 1\, \se h i$,\, $\se g i \in \se {{\cal
B}^1} m$\, and\, $\se h i \in \se {{\cal B}^1} k$, item (7.iii) yields\, $F(\se g i) = 
\se f m(\se g i)\, \spezz 2 \, \se f k(\se h i) = F(\se h i)$\;
($i = 1, \dots, r$). Then, by \ref{moreprods}\,(a) again,
\vspace{-0.15cm}

\nhp \hfl $F(g) = F(\se g 1) \cdot\ \dots\ \cdot F(\se g r)\; \spezz 2 \,
F(\se h 1) \cdot\ \dots\ \cdot F(\se h r) = F(h)$, \hfl
\vspace{-0.15cm}

\nhp which proves (c). The Claim follows from (a)--(c) using Lemma
\ref{fan-maps}\,(3). This completes the proof of Theorem \ref{isothffans}.
\hfl $\Box$

\vspace{-0.2cm}

\nhp \parbox[t]{220pt}{M. Dickmann\\
Institut de Math\'ematiques de Jussieu --\\
\hspace*{3.3cm} Paris Rive Gauche\\
Universit\'es Paris 6 et 7\\ 
Paris, France\\
e-mail: dickmann@math.univ-paris-diderot.fr} 
\nhp \hfl \parbox[t]{200pt}{A. Petrovich\\
Departamento de Matem\'atica\\
Facultad de Ciencias Exactas y Naturales\\
Universidad de Buenos Aires\\ 
Buenos Aires, Argentina\\
e-mail: apetrov@dm.uba.ar}

\end{document}